\input amstex
\input xypic
\documentstyle{amsppt}
\pagewidth{13.5cm}
\magnification=1200
\pageheight{20cm}
\define\alp{\alpha}
\define\be{\beta}
\def \jeden {1\hskip-3.5pt1}

\define\la{\lambda}
\define\La{\Lambda}

\define\ga{\gamma}
\define\Ga{\Gamma}
\define\de{\delta}

\define\bF{\Bbb F}
\define\bH{\Bbb H}

\define \bB{\Bbb B}
\define \bR{\Bbb R}
\define\bP{\Bbb P}
\define\bL{\Bbb L}
\define \bZ{\Bbb Z}
\define \bN{\Bbb N}
\define \bC{\Bbb C}

\define \bQ{\Bbb Q}
\def\op{\operatorname}

\define\gagin {\ga^{\op {Gin}}}
\topmatter
\title  Chern classes of proalgebraic varieties \\
and motivic measures
\endtitle
\author SHOJI YOKURA$^{*}$
\endauthor
\thanks {(*) Partially supported by Grant-in-Aid for Scientific Research (C) (No.15540086), the Japanese Ministry of Education, Science, Sports and Culture}  \newline \indent {\it Keywords:}  Bivariant Theory; Chern--Schwartz--MacPherson class; Constructible function; Grothendieck ring;
Inductive limit; Pro-object; Projective limit; Motivic measure\endthanks
\address
Department of Mathematics and Computer Science, 
Faculty of Science, 
University of Kagoshima, 21-35 Korimoto 1-chome, Kagoshima 890-0065, Japan
\endaddress
\email yokura$\@$sci.kagoshima-u.ac.jp
\endemail
\abstract {Michael Gromov has recently initiated what he calls ``symbolic algebraic geometry", in which objects are proalgebraic varieties: a proalgebraic variety is by definition the projective limit of a projective system of algebraic varieties. In this paper we construct Chern--Schwartz--MacPherson classes of proalgebraic varieties, by introducing the notion of ``proconstructible functions " and ``$\chi$-stable proconstructible functions" and using the Fulton-MacPherson's Bivariant Theory. As a ``motivic" version of a $\chi$-stable proconstructible function, $\Ga$-stable constructible functions are introduced. This construction naturally generalizes the so-called motivic measure and motivic integration. For the Nash arc space $\Cal L(X)$ of an algebraic variety $X$, the proconstructible set is equivalent to the so-called cylinder set or constructible set in the arc space.}
\endabstract 
\endtopmatter
\rightheadtext{Chern classes of proalgebraic varieties and motivic measures}
\document
\head \S 1 Introduction \endhead

This work is motivated by Gromov's papers [Grom 1, Grom 2] and also by [Y4]. 

A pro-category was first introduced by A. Grothendieck [Grot] and
 it was used to develope the Etale Homotopy Theory [AM] and Shape Theory (e.g., see [Bor], [Ed], [MS], etc.) and so on. 
A pro-algebraic variety is defined to be a projective system of complex algebraic varieties and 
a proalgebraic variety is defined to be the projective limit of a pro-algebraic variety. 
In [Grom 1] M. Gromov investigated the {\it surjunctivity} [Got] , {\it i.e., being either surjective or non-injective,}
 in the category of proalgebraic varieties. The original or classical surjunctivity theorem
is the so-called {\it Ax' Theorem}, saying that every regular selfmapping of a complex algebraic variety
 is surjunctive; thus if it is injective then it has to be surjective (cf, [Ax], [BBR], [Bo], [Kurd] , [New],  [Par], etc.).
 In [Grom 1] he initiated what he calles ``symbolic algebraic geometry" and in its Abstract he says
 `` ... The paper intends to bring out relations between model theory, algebraic geometry, and symbolic dynamics."

Our interest at the moment is not a further investigation concerning Ax-type theorems, but
characteristic classes, in particular, Chern classes of proalgebraic varieties. One motivation of this is that in [Grom 1] 
he uses the term of `proconstructible set or space' at several places; so we ask ourselves :
``whatever the definition of `proconstructible set' is, what would be
 the Chern--Schwartz--MacPherson class of a proalgebraic variety ?" This way of thinking is simply motivated by the fact
 that `constructible set' is nothing but `constructible function' and that in turn it is nothing but `homology class' 
via the Chern--Schwartz--MacPherson class transformation $c_*: F \to H_*$ (see [M], [BS], [Sc1, Sc2]). 
This transformation is the unique natural transformation from the covariant functor $F$ of constructible functions to 
the covariant homology functor $H_*$, satisfying the normalization condition that for a smooth variety $X$ 
the value $c_*(\jeden_X)$ of the characteristic function $\jeden_X$ is the Poincar\'e dual of
 the total Chern cohomology class $c^*(X)$. The unique existence of such a transformation was conjectured by Pierre Deligne and Alexander Grothendieck, and  finally was solved affirmatively by Robert MacPherson [M]. Later it was shown by Jean-Paul Brasselet and Marie-H\'el\`ene Schwartz [BS] that it is isomorphic to the Schwartz class [Sc1, Sc2] via the Alexander duality isomorphism. This is the reason for naming ``Chern--Schwartz--MacPherson".

Another motivation is that in [Y4] we investigated some ``quasi-bivariant" Chern classes, using the projective system of resolutions of singularities and its limits. (For bivariant theories, see [FM], [Br], [BSY], [E2], [EY1, 2], [F1], [G1], [Sa2], [Sch2], [Y1, 2, 3], [Z1, 2], etc.). What we did in [Y4] can be put in as follows: we make some bivariant classes obtained by all the different resolutions of singularities ``equivalent" or ``the same" in some inductive limit. 

A very simple example of a proalgebraic variety is the Cartesian product $X^{\bN}$ of an infinite countable copies of a complex algebraic variety $X$, which is one of the main objects treated in [Grom 1], i.e., a proalgebriac variety associated with a symbolic dynamic. Then, what would be
 the Chern--Schwartz--MacPherson of $X^{\bN}$ ?
In particular, what would be the {\it ``Euler-Poincar\'e (pro)characteristic"} of $X^{\bN}$ ?
Our answers are that they are respectively $c_*(X)$ and $\chi(X)$ in some sense, which will be clarified later.  It is this very simple observation (which looked stupid, nonsensical or meaningless at the beginning) that led us to the present work, which naturally led us to motivic measures, which have been actively studied by many people(e.g., see [Cr], [DL 1], [DL 2], [Kon], [Loo], [Ve] etc.).

Taking the degree of the 0-dimensional component or taking integration of the Chern--Schwartz--MacPherson class homomorphism $c_*: F(X) \to H_*(X)$ gives us the Euler Poincar\'e characteristic (homomorphism)
$$\chi : F(X) \to \bZ \quad \text {described by} \quad \chi(\alp) = \sum _{n \in \bZ} n \chi \left (\alp^{-1}(n) \right ). \tag 1.1$$
Another notation for this, putting emphasis on integration, is
$$ \int_X \alp d\chi = \sum _{n \in \bZ} n \chi \left (\alp^{-1}(n) \right ). \tag 1.1.a$$
Also, for a function $f(\alp)$, we can define the following
$$ \int_X f(\alp) d\chi = \sum _{n \in \bZ} f(n) \chi \left (\alp^{-1}(n) \right ). \tag 1.1.b$$
Which is a generalized version of (1.1.a).

Since a constructible function $\alp$ on $X$ can be expressed as a finite linear combination of characteristic functions $\jeden_W$ of subvarieties $W$ of $X$, $\alp = \sum n_W \jeden_W$, (1.1) is simply expressed as $\chi \left (\sum n_W \jeden_W \right ) = \sum n_W \chi(W).$ Another distinguished component of the $c_*$ is the ``top"-dimensional component. For $\alp = \sum n_W \jeden_W$ we have $c_*(\alp) = \sum n_W c_*(W).$ Since the top-dimensional component of each $c_*(W)$ is the fundamental class $[W]$, the ``top"-dimensional part or the ``fundamental class" part of $c_*(\alp)$ shall be denoted by $[c_*(\alp)]$ and $[c_*(\alp)] = \sum n_W [W].$
So, if we ``regard" the fundamental class $[W]$ as a class $[W]$ of the variety $W$ in the Grothendieck ring $K_0(\Cal V_{\bC})$ of complex algebraic varieties, then we can ``get" the homomorphism $\Ga :F(X) \overset {c_*} \to \longrightarrow H_*(X) \overset {[\quad]}  \to \longrightarrow K_0(\Cal V_{\bC})$.
This way of thinking is, however, not quite right, because the homomorphism $[\quad]: H_*(X) \to K_0(\Cal V_{\bC})$ is not well-defined. But, the homomorphism
$$\Ga :F(X) \to K_0(\Cal V_{\bC}) \quad 
\text {defined by} \quad 
\Ga (\alp) = \sum _{n \in \bZ} n \left [\alp^{-1}(n) \right ] \tag 1.2 $$
is a well-defined homomorphism, which is called the {\it Grothendieck class homomorphism}. Simply expressing, it is $\Ga \left (\sum n_W \jeden_W \right ) = \sum n_W [W].$ Thus (1.2) is a ``motivic" version of (1.1). Or as a ``motivic" version of (1.1.a), (1.2) can be expressed by $\displaystyle \int_X \alp d\Ga= \sum _{n \in \bZ} n \left [\alp^{-1}(n) \right ].$
And like (1.1.b), for a function $f(\alp)$ we have
$$\int_X f(\alp) d\Ga= \sum _{n \in \bZ} f(n) \left [\alp^{-1}(n) \right ].$$

Mimicking or abusing the above notations, we can also consider
$\displaystyle c_*(\alp) = \int_X \alp dc_* = \sum_{n \in \bZ} n c_* \left (\alp^{-1}(n) \right ).$
If we let $c_i: F \to H_{2i}$ be the $2i$-dimensional component of $c_*:F \to H_*$, we can also have $\displaystyle c_i(\alp) = \int_X \alp dc_i = \sum_{n \in \bZ} n c_i \left (\alp^{-1}(n) \right ).$

What we do in this paper is to extend or generalize characteristic classes, in particular, the Chern--Schwartz--MacPherson class $c_*: F(X) \to H_*(X)$, of complex (possibly singular) algebraic varieties to the category of proalgebraic varieties. The first thing that we have to consider is how to define the proalgebraic version of $F(X)$, namely, how to define a reasonable notion of {\it ``proconstructible function"}, i.e., a proalgebraic version of a constructible function. 

Especially, one can extend the above two homomorphisms $\chi: F(X) \to \bZ$ and $\Ga: F(X) \to K_0(\Cal V_{\bC})$ to the category of proalgebraic varieties and it turns out that the extended version include the so-called motivic measures as special cases. The key for extending these is that they are both multiplicative, while the other individual components of $c_*$ are not multiplicative and thus one cannot extend those to the category of proalgebraic varieties. 

The organization of the paper is as follows. In \S2 we just consider Chern class of a projective system of algebraic varieties, which is very straighforward. In \S 3 we define {\it `` proconstructible functions"} and we formulate Chern--Schwartz--MacPherson classes of some proalgebraic varieties, making a full use of the Bivariant Theory introduced by William Fulton and Robert MacPherson [FM], in particular a bivariant Chern class [Br]. 
The key of the results obtained in \S 3 is the constancy of the Euler--Poincar\'e characteristics of 
the fibers of each structure morphism in a projective system of varieties. 
If this constancy is not satisfied, we need a stronger requirement on the proconstructible function, and thus in \S 4 we introduce the notion of {\it $\chi$-stable proconstructible functions}. 
In \S 5 we discuss motivic measures and furthemore, hinted by the definition of $\chi$-stable proconstructible functions, we introduce the notion of 
{\it $\Ga$-stable (or ``motivic" stable) proconstructible functions}, and
thus we get a generalization of so-called {\it stable constructible functions}.  In \S6 we discuss the proalgebriac version of integration with respect to the proalgebraic Euler-Poincar\'e characteristic $\chi^{\op {pro}}$ and the proalgebraic Grothendieck ``motivic" class homomorphism $\Ga^{\op {pro}}$. In \S 7 we introduce the notion of {\it proresolution of singularities} and 
we pose a problem concerning a relationship between the Nash arc space $\Cal L(X)$ and
 the proresolution of $X$, which could be related to the problem which Nash considered 
in his paper [Na], in relation with Hironaka's resolution of singularities [Hi].

Finally we point out that a natural question coming out of this way of thinking is whether
one can obtain a theory of Chern--Schwartz--MacPherson classes with values in the Grothendieck ring, i.e., a ``motivic" version  of the Chern--Schwartz--MacPherson class $c_*$, and furthermore, if we consider the category of proalgebraic varieties, a natural question is whether
one can obtain a theory of Chern--Schwartz--MacPherson classes with values in localizations of the Grothendieck ring. These remain to be seen. And for another connection of Chern-Schwartz-MacPherson class with motivic measure or integration, see Paolo Aluffi's recent article [A].

{\it Acknowledgement.} The author would like to thank J\"org Sch\"urmann and Willem Veys for their many valuable comments and suggestions.

\head \S 2 Chern classes of pro-algebraic varieties \endhead

Let $I$ be a directed set and let $\Cal C$ be a given category. 
Then a projective system is, by definition, a system $\{X_i, \pi_{ii'}:X_{i'} \to X_i (i < i'), I \}$ 
consisting of objects $X_i \in \op {Obj}(\Cal C)$, morphisms 
$\pi_{ii'}:X_{i'} \to X_i  \in \op {Mor}(\Cal C)$ for each $i < i'$ and the index set $I$. 
The object $X_i$ is called a {\it term} and the morphism  $\pi_{ii'}:X_{i'} \to X_i$ 
a {\it bonding morphism} or {\it structure morphism} ([MS]). 
The projective system  $\{X_i, \pi_{ii'}:X_{i'} \to X_i (i < i'), I \}$ is
 sometimes simply denoted by $\{X_i\}_{i \in I}$.

Given a category $\Cal C$, Pro-$\Cal C$ is the category whose objects are projective systems 
$X = \{X_i\}_{i \in I}$ in $\Cal C$ and whose set of morphisms from 
$X = \{X_i \}_{i \in I}$ to $Y = \{Y_j \}_{j \in J}$ is

$$ \op {Pro-}\Cal C (X, Y) := \varprojlim_J (\varinjlim_I \Cal C(X_i, Y_j)).$$

This definition is not crystal clear, but a more down-to-earth definition is the following (e.g., see [Fox] or [MS]):
 A morphism $f:X \to Y$ consists of a map $\theta: J \to I$ (not necessarily order preserving) and morphisms
 $f_j: X_{\theta (j)} \to Y_j$ for each $j \in J$, subject to the condition that if $j < j'$ in 
$J$ then for some $i \in I$ such that $i> \theta (j)$ and $i > \theta (j')$, 
the following diagram commutes 

$$
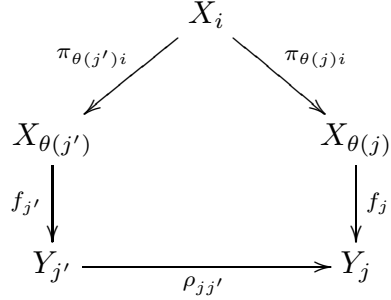
\xymatrix{
& X_i \ar[dl]_{\pi_{\theta (j') i}} \ar[dr]^ {\pi_{\theta (j)i}} & \\
X_{\theta(j')}\ar[d]_{f_{j'}} & & X_{\theta (j)} \ar[d]^{f_j} \\
Y_{j'} \ar[rr]_{\rho_{jj'}} && Y_j
}$$

Given a projective system $X = \{X_i \}_{i \in I} \in \op {Pro-} \Cal C$, 
the projective limit $X_{\infty} := \varprojlim X_i$ may not belong to the source category $\Cal C$. 
For a certain sufficient condition for the existence of the projective limit in the category $\Cal C$, 
see [MS] for example. 

An object in $\op {Pro-}\Cal C$ is called a {\it pro-object}. A projective system of algebraic varieties is called a {\it pro-algebraic
variety} and its projective limit is called a {\it proalgebraic variety}, which may not be an algebraic variety but simply a topological space.

Let $T: \Cal C \to \Cal D$ be a covariant functor between two categories $\Cal C, \Cal D$. Obviously the covariant functor $T$ extends to a covariant pro-functor
$$\op {Pro-}T : \op {Pro-}\Cal C \to \op {Pro-}\Cal D$$
defined by $\op {Pro-}T(\{X_i \}_{i \in I}) := \{T(X_i) \}_{i \in I}$. Let $T_1, T_2: \Cal C \to \Cal D$ be two covariant functors and $N: T_1 \to T_2$ be a natural transformation between the two functors $T_1$ and $T_2$. Then the natural transformation $N: T_1 \to T_2$ extends to a natural pro-transformation
$$ \op {Pro-} N : \op {Pro-} T_1 \to \op {Pro-} T_2.$$

Thus a pro-algebraic version of the Chern--Schwartz--MacPherson class is straightforward, i.e., we have
$$ \op {Pro-} c_*:  \op {Pro-} F \to  \op {Pro-} H_*$$
and thus there is nothing to be done for it.  In this case, the characteristic pro-function $\jeden_X$ of the pro-algebraic variety $X = \{X_i \}_{i \in I}$ should be simply $\jeden_X :=
 \{\jeden_{X_i} \}_{i \in I}$ and thus the pro-version of the Chern--Schwartz--MacPherson class of the pro-algebraic variety $X = \{X_i \}_{i \in I}$ is simply $\op {Pro-}c_*(X) =  \{c_*(X_i) \}_{i \in I}$.

What we want to do is its proalgebraic version.

\remark {Remark (2.1)} In Etale Homotopy Theory [AM] and Shape Theory (e.g., see [Bor], [Ed], [MS]) they stay in the pro-category and do not consider limits and colimits, because doing so throw away some geometric informations.
\endremark

\head \S3 Proconstructible functions and Chern classes of proalgebraic varieties \endhead

As mentioned  above, a pro-morphism between two pro-objects is quite complicated. However, it follows from [MS] that the pro-morphism can be described more naturally as a so-called {\it level preserving pro-morphism}. Suppose that we have two pro-algebraic varieties
$X = \{X_{\ga} \}_{\ga \in \Ga}$ and  $Y = \{Y_{\la} \}_{\la \in \La}$. Then a pro-algebraic morphism $\Phi =\{f_{\la} \}_{\la \in \La}: X \to Y$ is described as follows: there is an order-preserving map $\xi: \La \to \Ga$, i.e., $\xi (\la) < \xi (\mu)$ for  $\la < \mu $, and for each $\la \in \La$ there is a morphism $f_{\la}:  X_{\xi(\la)} \to Y_{\la}$ such that for $\la < \mu$ the following diagram commutes:

$$\CD
X_{\xi (\mu)}@> {f_{\mu}} >> Y_{\mu}\\
@V {\pi_{\xi(\la) \xi(\mu)}}VV @VV {\rho_{\la \mu}}V\\
X_{\xi (\la)}@>> {f_{\la}} > Y_{\la}, \endCD
$$

Then, the projective limit of the system $\{f_{\la} \}$ is a morphism from the proalgebraic variety $X_{\infty} = \varprojlim _{\la \in \La} X_{\la}$  to the proalgebraic variety $Y_{\infty} = \varprojlim _{\ga \in \Ga} Y_{\ga}$. It is called a proalgebraic morphism and denoted by
$f_{\infty}:X_{\infty} \to Y_{\infty}$.

The projective system $\{X_{\ga}, \pi_{\ga \de} (\ga < \de) \}$ induces the projective system of
abelian groups of constructible functions $\Bigl \{F(X_{\ga}), (\pi_{\ga \de})_*: F(X_{\de}) \to F(X_{\ga}) (\ga < \de) \Bigr \}$.
And a system of morphims $f_{\la}:  X_{\xi(\la)} \to Y_{\la}$ induces the system of homomorphisms ${f_{\la}}_*:  F(X_{\xi(\la)}) \to F(Y_{\la})$.
Thus the system of commutative diagrams

$$\CD
F(X_{\xi (\mu)})@> {{f_{\mu}}_*} >> F(Y_{\mu})\\
@V {{\pi_{\xi(\la) \xi(\mu)}}_*}VV @VV {{\rho_{\la \mu}}_*}V\\
F(X_{\xi (\la)})@>> {{f_{\la}}_*} > F(Y_{\la}), \endCD
$$
induces  the homomorphism ${f_*}_{\infty} : \varprojlim_{\ga \in \Ga} F(X_{\ga}) \to \varprojlim_{\la \in \La} F(Y_{\la})$.
Similarly we get the homomorphism of the projective limits of homology groups
${f_*}_{\infty} : \varprojlim_{\ga \in \Ga} H_*(X_{\ga}) \to \varprojlim_{\la \in \La} H_*(Y_{\la})$.
Furthermore the commutative diagram of Chern--Schwartz--MacPherson class homomorphisms

$$\CD
F(X_{\mu})@> {c_*} >> H_*(X_{\mu})\\
@V {{\pi_{\la \mu}}_*}VV @VV {{\pi_{\la \mu}}_*}V\\
F(X_{\la})@>> {c_*} > H_*(X_{\la}), \endCD
$$
induces the projective limit of Chern--Schwartz--MacPherson classes:
$${c_*}_{\infty} : \varprojlim_{\la \in \La} F(X_{\la}) \to \varprojlim_{\la \in \La} H_*(X_{\la})$$
So, we define, for the proalgebraic variety $X_{\infty} = \varprojlim_{\la \in \La} X_{\la}$,
$$\op {pro} F(X_{\infty}) := \varprojlim_{\la \in \La} F(X_{\la}) \qquad \text {and} \qquad \op {pro}H_*(X_{\infty}) := \varprojlim_{\la \in \La} H_*(X_{\la}).$$
If we define $\op {pro}c_*: \op {pro}F \to \op {pro}H_*$ to be the above ${c_*}_{\infty}$ and define $f_{\infty *}$ to be the above $f_{* \infty}$, then we have a na\" \i ve proalgebraic version of the Chern--Schwartz--MacPherson class

$$\op {pro}c_*: \op {pro}F \to \op {pro}H_*,$$

i.e., for a proalgebraic morphism $f_{\infty} :X_{\infty} \to Y_{\infty}$ we have the commutative diagram

$$\CD
\op {pro}F(X_{\infty})@> {\op {pro}c_*} >> \op {pro}H_*(X_{\infty})\\
@V {{f_{\infty}}_*}VV @VV {{f_{\infty}}_*}V\\
\op {pro}F(Y_{\infty})@>> {\op {pro}c_*} > \op {pro}H_*(Y_{\infty}). \endCD
$$

Note that if the directed set $\La$ is finite, then it is clear that there is an element $0 \in \La$ which is cofinal. In this case, the proalgebraic variety $X_{\infty}$ is isomorphic to the variety $X_0$ and thus the above na\" \i ve proalgebraic Chern--Schwartz--MacPherson class $\op {pro}c_*: \op {pro}F \to \op {pro}H_*$ is nothing but the original one. Thus we are of course interested in the case where the directed set $\La$ is infinite. 
For example, let the directed set be the natural numbers $\bN$ and the projective system is the sequence of complex algebraic varieties:
$$ \cdots \to X_i \to \cdots \to X_3 \to X_2 \to X_1.$$
For instance, for an algebraic variety $X$ of finite dimension, let us denote the Cartesian product $X \times X \times \cdots \times X$ of $n$ copies of $X$ by $X^n$. From now on, the superscript $n$ of $X^n$ does not mean the dimension of the variety $X$, unless stated otherwise. And we let $\pi_{n-1, n}:X^n \to X^{n-1}$ be the canonical projection defined by $\pi_{n-1,n}(x_1, x_2, \cdots, x_n) = (x_1, x_2, \cdots, x_{n-1})$ and consider the infinite
sequence
$$ \cdots \to X^n \to X^{n-1} \to \cdots \to X^3 \to X^2 \to X.$$
Its projective limit is the infinite product $X^{\bN} = X^{\infty}$, which is an important object in Gromov's paper [Grom 1]. The progroup $\displaystyle \op {pro}F(X^{\bN}) =  \varprojlim_n F(X^n)$ is,
 by definition

$$\op {pro}F(X^{\bN}) = \left \{ (\alp_j) \in \prod_{j =1}^{\infty}F(X^j) \quad \vert \quad \forall j \quad \alp_j = ({\pi_{j,j+1}})_* \alp_{j+1} \right \}.$$

The j-th component $\alp \in F(X^j)$ of $(\alp_j) \in \prod_{j =1}^{\infty}F(X^j)$ automatically determines the lower components $\alp_1, \alp_2, \cdots, \alp_{j-1}$, i.e., for each $k <j$

$$\alp_k = (\pi_{k,k+1})_* \cdots (\pi_{j-2,j-1})_* (\pi_{j-1,j})_* \alp.$$ 

However, as to the upper components $\alp_k (k > j)$, there are infinitely many choices and there is no formula to determine it. Thus, the structure of the progroup $\op {pro} F(X^{\infty})$ is not so obvious.

Next, how do we capture an element of $\varprojlim_{\la \in \La} F(X_{\la})$ as a function on the proalgebraic variety $X_{\infty} = \varprojlim_{\la \in \La} X_{\la}$ , namely, what is the value of an element $\alp_{\infty} = (\alp_{\la}) \in \varprojlim_{\la \in \La} F(X_{\la}) \subset \prod F(X_{\la})$ at a point
$(x_{\la}) \in \varprojlim_{\la \in \La} X_{\la} \subset \prod X_{\la}$ ?  A very na\" \i ve definition could be 
$$ \alp_{\infty}((x_{\la})) = \alp_{\la}(x_{\la}),$$
which is, however, not well-defined. Indeed, by the definition of the projective limit we have that $\pi_{\la \mu}(x_{\mu}) = x_{\la}$ and that ${\pi_{\la \mu}}_*(\alp_{\mu}) = \alp_{\la}$ for $\la < \mu$. The equality ${\pi_{\la \mu}}_*(\alp_{\mu}) = \alp_{\la}$ implies that 
$$\alp_{\la} (x_{\la}) = {\pi_{\la \mu}}_*(\alp_{\mu})(x_{\la}) = \chi \Bigl (\pi_{\la \mu}^{-1}(x_{\la}); \alp_{\mu} \Bigr )$$
which is the weighted Euler--Poincar\'e characteristic of the fiber $\pi_{\la \mu}^{-1}(x_{\la})$, or the sum of the 0-dimensional component of the Chern--Schwartz--MacPherson class $c_* \bigl (\alp_{\mu}|_{\pi_{\la \mu}^{-1}(x_{\la})} \bigr )$ of the constructible function $\alp_{\mu}|_{\pi_{\la \mu}^{-1}(x_{\la})}$ on the fiber $\pi_{\la \mu}^{-1}(x_{\la})$. Here, $\chi(A;\alpha) := \sum _{n \in \bZ} n \chi(A \cap \alp^{-1}(n))$. Thus, in general, these two equalities $\pi_{\la \mu}(x_{\mu}) = x_{\la}$ and ${\pi_{\la \mu}}_*(\alp_{\mu}) = \alp_{\la}$ do not imply that $\alp_{\la}(x_{\la}) = \alp_{\mu}(x_{\mu})$. Thus $\alp_{\infty}((x_{\la})) = \alp_{\la}(x_{\la})$ is not well-defined. Hence, an element of the progroup $\op {pro} F(X_{\infty})$ defined above would not be a good candidate to be considered as a function on the proalgebraic variety $X_{\infty}$.

However, the equalities $\pi_{\la \mu}(x_{\mu}) = x_{\la}$ and $\pi_{\la \mu}^*(\alp_{\la}) = \alp_{\mu}$ imply that $\alp_{\infty}((x_{\la})) = \alp_{\la}(x_{\la})$ is well-defined, since we have that 
$$\alp_{\mu}(x_{\mu}) = (\pi_{\la \mu}^*(\alp_{\la})) (x_{\mu}) = \alp_{\la} (\pi_{\la \mu} (x_{\mu})) = \alp_{\la}(x_{\la}).$$

So, it is reasonable to define the following

\definition {Definition (3.1)} For a proalgebraic variety $\displaystyle X_{\infty} = \varprojlim _{\la \in \La} X_{\la}$, the inductive limit of the inductive system $\Bigl \{F(X_{\la}), \pi_{\la \mu}^*: F(X_{\la}) \to F(X_{\mu}) (\la < \mu) \Bigr \}$ is denoted by $F^{\op {pro}}(X_{\infty})$, i.e.,
$$F^{\op {pro}}(X_{\infty}) := \varinjlim _{\la \in \La} F(X_{\la}) = \bigcup_{\mu} \rho^{\mu} \bigl (F(X_{\mu}) \bigr )$$
where $\rho^{\mu}: F(X_{\mu}) \to \varinjlim _{\la \in \La} F(X_{\la})$ is the homomorphism sending $\alp_{\mu}$ to its equivalence class $[\alp_{\mu}]$ of $\alp_{\mu}$. An element of the group $F^{\op {pro}}(X_{\infty})$ is called a {\it proconstructible} function on the proalgebraic variety $X_{\infty}$.
\enddefinition

The proconstructible function $[\jeden_{X_{\la}}]$ for any $\la \in \La$ shall be called {\it the procharacteristic function} on $X_{\infty}$ and denoted by $\jeden_{X_{\infty}}$.

The terminology {\it proconstructible} is used in [Grom 1], but its definition does not seem to be given explicitly in his paper.

Definition (3.1) can be used for any contravariant functor on the category of objects. Namely, if 
$\Cal F: \Cal C \to \Cal C'$ is a contravariant functor and $\bigl \{X_{\la}, \pi_{\la \mu}: X _{\mu} \to X_{\la} (\la < \mu) \bigr \}$ is a projective system in $\Cal C$, then for the projective limit
$\displaystyle X_{\infty} = \varprojlim _{\la \in \La} X_{\la}$, which itself {\it may not belong to} the category $\Cal C$, we can define
$$\Cal F^{\op {pro}}(X_{\infty}) := \varinjlim _{\la \in \La} \Bigl \{\Cal F(X_{\la}), \pi_{\la \mu}^*: \Cal F(X_{\la}) \to \Cal F(X_{\mu}) (\la < \mu) \Bigr \},$$
which also {\it may not belong to} the category $\Cal C'$. Although  $\Cal F^{\op {pro}}(X_{\infty})$ is in general an abstract object assigned to the projective limit $X_{\infty}$, in the case of constructible functions, $F^{\op {pro}}(X_{\infty})$ can be treated as a group of functions on  $X_{\infty}$ as observed above and more detailed discussion on $F^{\op {pro}}(X_{\infty})$ will be given in connection with motivic measures in a later section. In general situations, this abstract object is sufficient.

Now, instead of considering $F^{\op {pro}}(X_{\infty})$ of an arbitrary proalgebraic variety $X_{\infty}$, we first consider it for the infinite countable product $X^{\bN}$ of a complex algebraic variety $X$ as a simple model case.

For each projection $\pi_{(n-1)n}:X^n \to X^{n-1}$ , the pullback homomorphism $\pi_{(n-1)n}^*: F(X^{n-1}) \to  F(X^n)$ is the multiplication by the characteristic function $\jeden_X$ of the last factor $X$, i.e., 
$$\pi_{(n-1)n}^*(\alp) = \alp \times \jeden_X,$$
where  $(\alp \times \jeden_X)(y, x) := \alp(y) \jeden_X(x) = \alp (y).$ Then, using the cross product formula $c_*(\de \times \omega) = c_*(\de) \times c_*(\omega)$ of the Chern--Schwartz--MacPherson class $c_*$ (see [Kw] and also cf. [KY]), we get the following commutative diagram

$$
\CD 
F(X^{n-1}) @ > {\times \jeden_X}>> F(X^n) \\
@V {c_*} VV  @ VV {c_*} V \\
H_*(X^{n-1}) @ >> {\times c_*(X)} > H_*(X^n). \endCD
$$

So, if we set
$$H_*^{\op {pro}}(X^{\bN}) := \varinjlim_n \Bigl \{\times c_*(X):H_*(X^{n-1}) \to H_*(X^n) \Bigr \} ,$$
then we have a {\it proalgebraic Chern--Schwartz--MacPherson class homomorphism}:
$$c_*^{\op {pro}}: F^{\op {pro}}(X^{\bN})  \to H_*^{\op {pro}}(X^{\bN})$$
Let us look at the ``0-dimensional component" of this homomorphism, i.e., 
$\displaystyle \int_{X^{\bN}} c_*^{\op {pro}}$. Namely, we consider the following commutative diagram for each$n$

$$
\CD 
F(X^{n-1}) @ > {\times \jeden_X} >> F(X^n) \\
@V {\chi} VV  @ VV {\chi} V \\
\bZ @>> {\times \chi(X)} > \bZ. \endCD \tag 3.2
$$
and we take its inductive limit, which gives us {\it the proalgebraic Euler--Poincar\'e characteristic (homomorphism)} 
$$\chi^{\op {pro}}: F^{\op {pro}}(X^{\bN})  \to \varinjlim_n \Bigl \{\times \chi(X): \bZ \to \bZ \Bigr \}.$$

It is clear from the definition that for $\alp_n \in F(X^n)$ we have $\chi^{\op {pro}}([\alp_n]) = [\chi(\alp_n)]$. To describe $[\chi(\alp_n)]$, i.e., $\chi^{\op {pro}}$,  more explicitly or down-to-earth, we recall the following lemma:
\proclaim {Lemma (3.3)} Let $p$ be a non-zero integer. For each positive integer $n$, let $G_n = \bZ$ be the integers and $\pi_{n(n+1)}: G_n \to G_{n+1}$ be the homomorphism defined by multiplication by $p$, i.e., $\pi_{n(n+1)}(m) = pm$. Then the inductive limit $\varinjlim _n G_n$ is isomorphic to the abelian group of all rationals of the form $\displaystyle \frac {m}{p^n}$, for integers $m$ and $n$, which shall be denoted by $\displaystyle \bZ \left [\frac {1}{p} \right ]$. If $p = 0$, then the inductive limit  $\varinjlim _n G_n = 0.$ \qed
\endproclaim

This seems to be a well-known fact (cf. [Munk, \S 73, Example 4, p. 435]) and this simple fact turns out to be a key in what follows. For the sake of later use, we give a quick proof. For each $n$, let
$\phi^n: G_n \to \bZ \left [\frac {1}{p} \right ]$ be defined by $\displaystyle \phi^n(m) = \frac {m}{p^{n-1}}$. Then we have the commutative diagram

$$\xymatrix{
G_n \ar[dr]_ {\phi^n }\ar[rr]^ {\pi_{n(n+1)}} && G_{n+1}\ar[dl]^{\phi^{n+1}}\\
& \bZ \left [\frac {1}{p} \right ]
}$$

It follows from standard facts on inductive systems and inductive limits that all these commutative diagrams give rise to the unique homomorphism 
$$\Phi: \varinjlim_n \phi^n : \varinjlim _n G_n \to \bZ \left [\frac {1}{p} \right ] \quad \text {described by} \quad 
\Phi \left (\sum_k \rho^k(m_k) \right ) = \sum_k \frac {m_k}{p^{k-1}},$$
which turns out to be surjective and injective.

There are of course infinitely many other isomorphims between $\varinjlim _n G_n$ and $\bZ \left [\frac {1}{p} \right ]$. Indeed, let for any integer $w$, we define $\phi^n_w: G_n \to \bZ \left [\frac {1}{p} \right ]$ by $\displaystyle \phi^n_w(m) = \frac {m}{p^{n-1+w}}$. Then for the integer $w$ we get an isomorphism 

$$\Phi_w = \varinjlim_n \phi^n_w : \varinjlim _n G_n \to \bZ \left [\frac {1}{p} \right ] \quad \text {defined by} \quad  \Phi_w \left (\sum_k \rho^k(m_k) \right ) = \sum_k \frac {m_k}{p^{k-1+w}}.$$
When $w = 0$, the above $\Phi_0 = \Phi = \varinjlim_n \phi^n : \varinjlim _n G_n \to \bZ \left [\frac {1}{p} \right ]$ is called the {\it canonical isomorphism}.

Using this lemma we have the following proalgebraic version of the Euler--Poincar\'e characteristic homomorphism $\chi: F(V) \to \bZ$:

\proclaim {Theorem (3.4)} Let $X$ be a complex algebraic variety and we assume that $\chi(X) \not = 0$. The canonical Euler--Poincar\'e (pro)characteristic homomorphism
$$\chi^{\op {pro}}: F^{\op {pro}}(X^{\bN}) \to \bZ \left [\frac {1}{\chi (X)} \right ]$$
is described by
$$\chi^{\op {pro}}\left (\sum_n [\alp_n] \right ) = \sum_n \frac {\chi(\alp_n)}{\chi(X)^{n-1}} $$
where $\alp_n \in F(X^{n})$. \qed
\endproclaim

In particular, if $\jeden_{X^{\bN}}$ denotes the characteristic function on the proalgebraic variety $X^{\bN}$, hence $\jeden_{X^{\bN}}= [\jeden_{X^k}]$ for any $k$, we have

$$\chi^{\op {pro}}\left (\jeden_{X^{\bN}}\right ) = \frac{\chi (\jeden_{X^k})}{\chi(X)^{k -1}}= \frac{\chi(X)^k}{\chi(X)^{k -1}} = \chi (X)$$
which is called {\it the canonical Euler--Poincar\'e (pro)characteristic} of the proalgebraic variety $X^{\bN}$ and simply denoted by $\chi^{\op {pro}} (X^{\bN})$.

This simple and na\" \i ve observation, which looked nonsensical or meaningless at the beginning, was the very start of the present work. Note that $c_*^{\op {pro}}(X^{\bN}) := c_*^{\op {pro}}(\jeden_{X^{\bN}}) = [c_*(X)]$, which is of course equal to $[c_*(X^k)]$ for any positive integer $k$.

As pointed out above, for any integer $w$
$$\chi_w^{\op {pro}}\left (\sum_n [\alp_n] \right ) := \sum_n \frac {\chi(\alp_n)}{\chi(X)^{n -1+w}}  $$
is also an Euler--Poincar\'e (pro)characteristic homomorphism. Unless stated otherwise, we consider only the above canonical one.

Theorem (3.4) can be generalized to the following:

\proclaim {Theorem (3.5)} Let $\displaystyle X_{\infty}= \varprojlim_{n \in \bN}X_n$ be a proalgebraic variety such that for each $n$ the structure morphism $\pi_{n(n+1)}: X_{n+1} \to X_n$ satisfies the condition that the Euler--Poincar\'e characteristics of the fibers of $\pi_{n, n+1}$ are non-zero (which implies the surjectivity of the morphism $\pi_{n(n+1)}$) and the same; for example,  $\pi_{n(n+1)}: X_{n+1} \to X_n$ is a locally trivial fiber bundle with fiber variety being $F_n$ and $\chi (F_n) \not = 0$. Let us denote the constant Euler--Poincar\'e characteristic of the fibers of the morphism $\pi_{n(n+1)}: X_{n+1} \to X_n$ by $\chi_n$ and we set $\chi_0 :=1$.

(i) The canonical Euler--Poincar\'e (pro)characteristic homomorphism
$$\chi^{\op {pro}}: F^{\op {pro}}(X_{\infty}) \to \bQ$$
is described by
$$\chi^{\op {pro}}\left (\sum_n [\alp_n] \right ) = \sum_n \frac {\chi(\alp_n)}{\chi_0 \cdot \chi_1\cdot \chi_2 \cdots \chi_{n-1}}. $$

(ii) In particular, if the Euler-Poincar\'e characteristics $\chi_n$ are all the same, say $\chi_n = \chi$ for any $n$, then the canonical Euler--Poincar\'e (pro)characteristic homomorphism
$\chi^{\op {pro}}: F^{\op {pro}}(X_{\infty}) \to \bQ$
is described by
$$\chi^{\op {pro}}\left (\sum_n [\alp_n] \right ) = \sum_n \frac {\chi(\alp_n)}{\chi^{n-1}}. $$
In this special case, the target ring $\bQ$ can be replaced by the ring $\bZ \left  [\frac {1}{\chi} \right ]$, the non-canonical one $\chi_w^{\op {pro}}$ is described by
$$\chi_w^{\op {pro}}\left (\sum_n [\alp_n] \right ) = \sum_n \frac {\chi(\alp_n)}{\chi^{n-1+w}}. $$

(iii) In particular, if $\chi_n = \chi (\bP^n)$ for each $n$ or $F_n = \bP^n$ the complex projective space of dimension $n$ for each $n$ in the above example, then the canonical Euler--Poincar\'e (pro)characteristic homomorphism $\chi^{\op {pro}}$ is surjective.
 \qed
\endproclaim

In particular, by considering the characteristic function $\jeden_{X_{\infty}}= [\jeden_{X_n}]$ for any $n$,
$$\chi^{\op {pro}}\left (X_{\infty}\right ) := \chi^{\op {pro}}\left (\jeden_{X_{\infty}}\right )= \frac{\chi(X_1)\cdot \chi_1\cdot \chi_2\cdots \chi_{n-1}}{\chi_1\cdot \chi_2 \cdots \chi_{n-1}} = \chi (X_1),$$
which is called {\it the canonical Euler--Poincar\'e (pro)characteristic} of the above proalgebraic variety $X_{\infty}$.

\demo {Proof of Theorem (3.5)}For a morphism $f:X \to Y$ whose fibers all have the same non-zero Euler--Poincar\'e (pro)characteristic, denoted by $\chi_f$, we get the commutative diagram

$$
\CD 
F(Y) @ > {f^*} >> F(X) \\
@V {\chi} VV  @ VV {\chi} V \\
\bZ @>> {\times \chi_f} > \bZ. \endCD \tag 3.5.1
$$
To see the commutativity of this diagram, we use the fact that for any constructible function $\be \in F(X)$ we have
$$\chi(\be) = \chi (f_*\be),$$
which follows from the naturality of the Chern--Schwartz--MacPherson class $c_*$. 
Indeed, using the above fact, we can show that for any characteristic function
$\jeden_W \in F(Y)$ we have
$$
\align
\chi \bigl (f^*\jeden_W \bigr ) & = \chi \bigl (f_*f^*\jeden_W \bigr ) \\
 & = \chi (\chi_f \cdot \jeden_W) \\
& = \chi_f \cdot \chi (\jeden_W). \endalign
$$
Therefore the above diagram (3.5.1) commutes. Thus the theorem follows from the following commutative diagram for each $n$;

$$
\CD 
F(X_n) @ > {\pi_{n(n+1)}^*} >> F(X_{n+1}) \\
@V {\chi} VV  @ VV {\chi} V \\
\bZ @>> {\times \chi_n} > \bZ \endCD \tag 3.5.2
$$
and the following generalized version of Lemma (3.3):

\proclaim {Lemma (3.5.3)} For each positive integer $n$, let $G_n = \bZ$ be the integers and $\pi_{n, n+1}: G_n \to G_{n+1}$ be the homomorphism defined by multiplication by a non-zero integer $p_n$, i.e., $\pi_{n, n+1}(m) = mp_n$. Then there exists a unique (injective) homomorphism
$$\Psi:\varinjlim _n G_n \to \bQ$$
such that the following diagram commutes

$$\xymatrix{
& G_n \ar [dl]_{\rho^n} \ar [dr]^{\times \frac {1}{p_0p_1\cdots p_{n-1}}} \\
{\varinjlim _n G_n} \ar [rr] _{\Psi}& &  \Bbb Q.}
$$

Here we set $p_0:=1$. And $\Psi$ is described by
$$\Psi \left (\sum_n \rho^n(r_n) \right ) = \sum_n \frac {r_n}{p_0 p_1 \cdots p_{n-1}}. $$
\endproclaim
Thus we get (i) and (ii). In the case of (iii), the injective homomorphism $\Psi:\varinjlim _n G_n \to \bQ$ becomes surjective, thus we get (iii).
\qed
\enddemo

Let the situation be as in Theorem (3.5). Let $\displaystyle \widehat {F^{\op {pro}}}(X_{\infty})$ denote the abelian group of formal infinite sums $\displaystyle \sum_{n=1}^{\infty} [\alp_n]$ such that its proalgebraic Euler--Poincar\'e characteristic $\displaystyle \sum_{n=1}^{\infty} \chi^{\op {pro}} \bigl ([\alp_n] \bigr )$ converges. Thus $\displaystyle \widehat {F^{\op {pro}}}(X_{\infty})$ can be interpreted as a {\it completion} of $\displaystyle F^{\op {pro}}(X_{\infty})$. With this definition we can get the following theorem, which is impossible in the usual algebraic geometry.

\proclaim {Theorem (3.6)}For the infinite product space 
$$\prod_{n=1}^{\infty}\bP^n := \bP^1 \times \bP^2 \times \cdots \times \bP^n \times \cdots$$ 
(which is a proalgebraic variety), the proalgebraic Euler--Poincar\'e characteristic homomorphism
$$\widehat {\chi^{\op {pro}}}: \widehat  {F^{\op {pro}}} \left (\prod_{n=1}^{\infty}\bP^n \right ) \to \bR \quad \text {defined by} \quad \widehat {\chi^{\op {pro}}} \left (\sum_{n=1} ^{\infty}[\alp_n] \right ) = \sum_{n=1} ^{\infty}\chi^{\op {pro}} \bigl  ([\alp_n] \bigr )$$
is surjective.
\qed
\endproclaim

Up to now, to get our results, we use the commutative diagrams (3.2) and (3.5.2), i.e., the constancy of the Euler--Poincar\'e characteristics of all the fibers of each structure morphsim. Thus, in order to consider the proalgebraic Chern--Schwartz--MacPherson classes or the Euler--Poincar\'e characteristics of proalgebraic varieties from the viewpoint of inductive limits, structure morphisms constituting projective systems must have such a strong requirement.

In fact, Theorem (3.5) can be extended to $c_*^{\op {pro}}$. To state and prove such a theorem, we need to appeal to the Bivariant Theory introduced by William Fulton and Robert MacPherson [FM], in particular a bivariant Chern class [Br]. So, we quickly recall only necessary ingredients of the Bivariant Theory for our use.

A bivariant theory $\bB$ on a category $\Cal C$ with values in the category of abelian groups is an assignment to each morphism
$$ X  \overset f \to \longrightarrow Y$$
in the category $\Cal C$ a graded abelian group
$$\bB(X  \overset f \to \longrightarrow Y)$$
which is equipped with the following three basic operations:

\noindent (Product operations): For morphisms $f: X \to Y$ and $g: Y
\to Z$, the product operation
$$\bullet: \bB( X \overset f \to \longrightarrow Y) \otimes \bB( Y \overset g
\to \longrightarrow Z) \to
\bB( X \overset {gf} \to \longrightarrow Z)$$
is  defined.

\noindent (Pushforward operations): For morphisms $f: X \to Y$
and $g: Y \to Z$ with $f$ proper, the pushforward operation
$$f_{\bigstar}: \bB( X \overset {gf} \to \longrightarrow Z) \to \bB( Y \overset g
\to \longrightarrow Z) $$
is  defined.

\noindent (Pullback operations): For a fiber square
$$\CD
X' @> g' >> X \\
@V f' VV @VV f V\\
Y' @>> g > Y, \endCD
$$
the pullback operation
$$g^{\bigstar} : \bB( X \overset  f \to \longrightarrow Y) \to \bB( X' \overset {f'}
\to \longrightarrow Y') $$
is  defined.

And these three operations are required to satisfy the seven compatibility axioms (see [FM, Part I, \S 2.2] for details). 

Let $\bB, \bB'$ be two bivariant theories on a category $\Cal C$. Then
a {\it Grothendieck transformation} from $\bB$ to $\bB'$
$$\ga : \bB \to \bB'$$
is a collection of homomorphisms
$$\bB(X \to Y) \to \bB'(X \to Y) $$
for a morphism $X \to Y$ in the category $\Cal C$, which preserves the above three basic operations: 

\hskip1cm (i) \quad $\ga (\alp \bullet_{\bB} \be) = \ga (\alp) \bullet _{\bB'} \ga (\be)$, 

\hskip1cm  (ii) \quad $\ga(f_{\bigstar}\alp) = f_{\bigstar} \ga (\alp)$, and 

\hskip1cm  (iii) \quad $\ga (g^{\bigstar} \alp) = g^{\bigstar} \ga (\alp)$.

A bivariant theory unifies both a covariant theory and a contravariant theory in the following sense:
$B_*(X):= \bB(X \to pt)$ and $B^*(X) := \bB(X \overset  {id} \to \longrightarrow X)$ become a covariant functor and a contravariant functor, respectively. And a Grothendieck transformation $\ga: \bB \to \bB'$ induces natural transformations $\ga_*: B_* \to B_*'$ and $\ga^*: B^* \to {B'}^*$.
Note also that if we have a Grothendieck transformation $\ga: \bB \to \bB'$, then via a bivariant class $b \in \bB( X \overset  f \to \longrightarrow Y)$ we get the commutative diagram
$$\CD
B_*(Y) @> {\ga_*}>> B'_*(Y) \\
@V {b \bullet} VV @VV {\ga(b) \bullet }V\\
B_*(X) @>> {\ga_*}> B'_*(X). \endCD
$$
This is called {\it the Verdier-type Riemann--Roch associated to the bivariant class $b$}.

Fulton--MacPherson's bivariant group $\bF(X \overset f \to \longrightarrow Y)$ of constructible functions
 consists of all the constructible functions on $X$ which satisfy the local Euler condition with respect to $f$. 
 Here a constructible function $\alp \in F(X)$ is said to satisfy the {\it local Euler condition with respect to $f$}
if for any point $x \in X$ and for any local embedding $(X, x) \to (\bold C^N,
0)$ the equality $\alp(x) = \chi \left (B_{\epsilon} \cap f^{-1}(z);\alp \right)$ holds,
where $B_{\epsilon}$ is a sufficiently small open ball of the origin $0$
with radius $\epsilon$ and $z$ is any point close to $f(x)$ (cf. [Br], [Sa2]). 
In particular, if $\jeden_f := \jeden_X$ belongs to the bivariant group $\bF(X \overset f \to \longrightarrow Y)$, then the morphism $f: X \to Y$ is called {\it an Euler morphism}.
And any constructible function in the bivariant group $\bF(X \overset f \to \longrightarrow Y)$
is called {\it a bivariant constructible function} to emphasize the bivariantness.

The three operations on $\bF$
are defined as follows:

(i) the product operation
$\bullet: \bF( X \overset f \to \longrightarrow Y) \otimes \bF( Y \overset g
\to \longrightarrow Z) \to
\bF( X \overset {gf} \to \longrightarrow Z)$
is  defined by 
$$\alp \bullet \be:= \alp \cdot f^*\be$$, 

(ii) the pushforward operation
$f_{\bigstar}: \bF( X \overset {gf} \to \longrightarrow Z) \to \bF( Y \overset g
\to \longrightarrow Z) $
is the usual pushforward $f_*$, i.e., 
$$f_{\bigstar}(\alp)(y):= \int c_*(\alp|_{f^{-1}}),$$

(iii) for a fiber square
$$\CD
X' @> g' >> X \\
@V f' VV @VV f V\\
Y' @>> g > Y, \endCD
$$
the pullback operation
$g^{\bigstar} : \bF( X \overset  f \to \longrightarrow Y) \to \bF( X' \overset {f'}
\to \longrightarrow Y') $
is the functional pullback ${g'}^*$, i.e.., 
$$g^{\bigstar}(\alp)(x'):= \alp (g'(x')).$$

Note that $\bF(X \overset {\op {id}_X} \to \longrightarrow X)$ consists of all locally constant functions and $\bF(X \to pt) = F(X)$. As a corollary of this observation, we have

\proclaim {Proposition (3.7)} For any bivariant constructible function $\alp \in \bF( X \overset  f \to \longrightarrow Y)$, the Euler--Poincar\'e characteristic $\displaystyle \chi \bigl (f^{-1}(y); \alp \bigr ) = \int c_* \bigl (\alp|_{f^{-1}(y)} \bigr )$ of $\alp$ restricted to each fiber $f^{-1}(y)$ is locally constant, i.e., constant along connected components of the base variety $Y$. In particular, if $f: X \to Y$ is an Euler morphism, then the Euler--Poincar\'e characteristic of the fibers are locally constant.
\endproclaim 

Note that locally trivial fiber bundles are Euler, but not vice versa.

Let $\bH$ be Fulton--MacPherson's bivariant homology theory, constructed 
from the cohomology theory. For a morphism $f: X \to Y$, choose a morphism $\phi: X \to \bR^n$ such that 
$\Phi:= (f, \phi): X \to Y \times \bR^n$ is a closed embedding. Then the $i$-th bivariant homology group 
$\bH^i(X \overset f \to \longrightarrow Y)$ is defined by

$$  \bH^i(X \overset f \to \longrightarrow Y) := H^{i+n}(Y \times \bR^n, Y \times \bR^n \setminus X_{\phi}),$$
where $X_{\phi}$ is defined to be the image of the morphism $\Phi = (f,\phi)$. The definition  is independent of the choice of $\phi$. Note that instead of taking the Euclidean space $\bR^n$ we can take a manifold $M$ so that $i:X \to M$ is a closed embedding and then consider the graph embedding $f \times i: X \to Y \times M$. See [FM, \S 3.1] for more details of $\bH$. In particular, note that if $Y$ is a point $pt$, $\bH(X \to pt)$ is isomorphic to the homology group $H_*(X)$ of the source variety $X$. 

W. Fulton and R. MacPherson conjectured or posed as a question the existence of a so-called bivariant Chern class and J.-P. Brasselet [Br] solved it:

\proclaim {Theorem (3.8)} (J.-P. Brasselet) On the category of complex analytic varieties and cellular morphisms, there exists a Grothendieck transformation
$$\ga^{\op {Br}} : \bF \to \bH$$
satisfying the normalization condition that $\ga^{\op {Br}}(\jeden_{\pi}) = c(TX) \cap [X]$  for $X$ smooth, where $\pi: X \to pt$ and $\jeden_{\pi} = \jeden_X$.
\endproclaim

Note that for a morphism $f: X \to pt$ from a variety $X$ to a point $pt$, $\ga^{\op {Br}}: \bF(X \to pt) \to \bH(X \to pt)$ is nothing but the original Chern--Schwartz--MacPherson class $c_*: F(X) \to H_*(X).$

\proclaim {Corollary (3.9)} (Verdier-type Riemann--Roch for Chern class) For a bivariant constructible function $\alp \in \bF( X \overset f \to \longrightarrow Y)$ we have the following commutative diagram:
$$\CD
F(Y) @> {c_*}>> H_*(Y) \\
@V {\alp \bullet_{\bF} = \alp \cdot f^*} VV @VV {\ga^{\op {Br}} (\alp) \bullet_{\bH} }V\\
F(X) @>> {c_*}> H_*(X). \endCD
$$
In particular, for an Euler morphism we have the following diagram:
$$\CD
F(Y) @> {c_*}>> H_*(Y) \\
@V {\jeden_f \bullet_{\bF} = f^*} VV @VV {\ga^{\op {Br}} (\jeden_f) \bullet_{\bH} }V\\
F(X) @>> {c_*}> H_*(X). \endCD
$$
(The homomorphism $\ga^{\op {Br}} (\jeden_f) \bullet_{\bH}$ shall be denoted by $f^{**}$.)
\endproclaim
 For a more generalized Verdier-type Riemann--Roch theorem for Chern--Schwartz--MacPherson class, see J\"org Sch\"urmann's recent article [Sch1].

\proclaim {Theorem (3.10)} Let $\displaystyle X_{\infty} = \varprojlim_{\la \in \La}X_{\la}$ be 
a proalgebraic variety such that for each $\la < \mu$ the structure morphism $\pi_{\la \mu}: X_{\mu} \to X_{\la}$ is an Euler proper morphism (hence surjective) of topologically connected algebraic varieties with the constant Euler-Poincar\'e characteristic $\chi_{\la \mu}$ of the fiber of the morphism $\pi_{\la \mu}$ being non-zero. Let $H_*^{\op {pro}}(X_{\infty})$ be the inductive limit of the inductive system $\displaystyle \Bigl  \{ \pi_{\la \mu}^{**}: H_*(X_{\la}) \to H_*(X_{\mu}) \Bigr \}$.

There exists a proalgebraic Chern--Schwartz--MacPherson class homomorphism
$$c_*^{\op {pro}}: F^{\op {pro}}(X_{\infty})  \to H_*^{\op {pro}}(X_{\infty})
\quad \text {defined by} \quad c_*^{\op {pro}}\left (\sum_{\la} [\alp_{\la}] \right ) = \sum_{\la} \rho^{\la}\bigl (c_*(\alp_{\la}) \bigr ). $$
In particular, the proalgebraic integration becomes the following
$$\chi^{\op {pro}}: F^{\op {pro}}(X_{\infty})  \to \varinjlim_{\La} \Bigl \{\times \chi_{\la \mu}: \bZ \to \bZ \Bigr \}.$$
(Here we do not know an explicit description of the inductive limit $\displaystyle \varinjlim_{\La} \bigl \{\times \chi_{\la \mu}: \bZ \to \bZ \bigr  \}$ like Lemma (3.5.3).)
\qed
\endproclaim
 In this general situation we have that
$$c_*^{\op {pro}}(X_{\infty}) = [c_*(X_{\la})] \quad \text {and} \quad \chi^{\op {pro}}(X_{\infty}) = [\chi(X_{\la})] $$
for any $\la \in \La$, but we cannot give a more explicit description like in the case when $\La = \bN$.

What we have done so far is the proalgebraic Chern--Schwartz--MacPherson class homomorphism, and our eventual problem is whether one can capture this homomorphism as {\it a natural transformation} as in the original Chern--Schwartz--MacPherson class.

First, as a trial, we consider a very simple promorphism
$$f: X^{\bN} \to Y^{\bN}$$
which is one of the important objects in Gromov's papers [Grom 1, Grom 2]. Namely, we consider the naturality of $c_*^{\op {pro}}: F^{\op {pro}}(X^{\bN})  \to H_*^{\op {pro}}(X^{\bN})$. Even in this simple case it turns out that we cannot expect a reasonable result.  For example, for a promorphism $f: X^{\bN} \to Y^{\bN}$ consider the the projective limit of
$f^n :=\undersetbrace \text {$n$ times } \to { f \times \cdots \times f}: X^n \to Y^n$ for a morphsim $f: X \to Y$. For more nontrivial examples involving subtle combinatoric natures of graphs, see [Grom 1]. Even in the above simplest case, one cannot get a reasonable solution even if we consider the naturality of the homomorphism $\chi^{\op {pro}}: F^{\op {pro}}(X^{\bN}) \to \bZ \left [\frac {1}{\chi (X)} \right ]$. Our first simple-minded and hasty answer was that if we just took the ``denominator changing" homomorphism $D_{X/Y}: \bZ \left [\frac {1}{\chi (X)} \right ] \to \bZ \left  [\frac {1}{\chi (Y)} \right ]$ defined by 
$$D_{X/Y}\left ( \sum_k \frac {a_k}{\chi(X)^k}\right ) := \sum_k \frac {a_k}{\chi(Y)^k}$$
then we would get the following commutative diagram, i.e., the naturality of $\chi^{\op {pro}}$:
$$
\CD 
F^{\op {pro}}(X^{\bN}) @ > {f_*} >> F^{\op {pro}}(Y^{\bN})\\
@V {\chi^{\op {pro}}} VV  @ VV {\chi^{\op {pro}}} V \\
\bZ \left  [\frac {1}{\chi (X)} \right ]  @>> {D_{X/Y}} > \bZ \left  [\frac {1}{\chi (Y)} \right ] . \endCD
$$
Which is ``because" of the following computation:
$$
\align
\chi^{\op {pro}}\left ( f_* \left (\sum_n [\alp_n] \right ) \right ) & = \sum_n \frac {\chi \bigl ( (f^n)_*\alp_n \bigr )}{\chi(Y)^{n -1}} \\
& = \sum_n \frac {\chi (\alp_n)}{\chi(Y)^{n -1}} \quad \text {(since $\chi \bigl ((f^n)_*\alp_n \bigr ) = \chi (\alp_n)$) } \\
& = D_{X/Y} \left (\sum_n \frac {\chi (\alp_n)}{\chi(X)^{n -1}}\right ) \\
& = D_{X/Y} \left (\chi^{\op {pro}}\left ( \sum_n [\alp_n] \right ) \right ). \endalign
$$
However, this ``denominator changing" homomorphism $D_{X/Y}: \bZ \left [\frac {1}{\chi (X)} \right ] \to \bZ \left [\frac {1}{\chi (Y)} \right ]$ is not well-defined, since it is not independent of the expression of the rational.

So, to get a reasonable result, we have to restrict ourselves to some special promorphisms.

 If the commutative diagram

$$\CD
Y_{\xi (\mu)}@> {f_{\mu}} >> X_{\mu}\\
@V {\rho_{\xi(\la) \xi(\mu)}}VV @VV {\pi_{\la \mu}}V\\
Y_{\xi (\la)}@>> {f_{\la}} > X_{\la}, \endCD
$$

is a fiber square, then we call the pro-morphism $\{f_{\la}:  Y_{\xi(\la)} \to X_{\la} \}$ a {\it fiber-square pro-morphism}, abusing words.

Let $\{f_{\la}:  Y_{\xi(\la)} \to X_{\la} \}$ be a fiber-square pro-morphism between two pro-algebraic varieties with Euler morphisms of possibly singular terms $X_{\la}'s$ and $Y_{\ga}'s$. Then our question is whether the following diagram is commutative or not:

$$\CD
F^{\op {pro}}(Y_{\infty}) @> {c_*^{\op {pro}}} >>H_*^{\op {pro}}(Y_{\infty})\\
@V {{f_{\infty}}_*}VV @VV {{f_{\infty}}_*}V\\
F^{\op {pro}}(X_{\infty}) @>> {c_*^{\op {pro}}} > H_*^{\op {pro}}(X_{\infty}). \endCD
$$

To prove the commutativity of this diagram, it suffices to show the commutativity of the following diagrams:

$$\diagram
& F(Y_{\xi (\la)}) \xto'[1,0][2,0]_(0.45){\rho_{\xi(\la) \xi(\mu)}^* }\dlto_{c_*} \rrto^{{f_{\la}}_*} && F(X_{\la})\ddto^{\pi_{\la \mu}^* }\dlto^{c_*} \\
H_*(Y_{\xi (\la)}) \ddto_{\rho_{\xi(\la) \xi(\mu)}^{**}} \rrto^(0.7){{f_{\la}}_*} && H_*(X_{\la}) \ \ddto^(0.35){\pi_{\la \mu}^{**}} \\
& F(Y_{\xi (\mu)}) \quad \dlto_{c_*} \xto'[0,1][0,2]_(-0.35){{f_{\mu}}_*}&& F(X_{\mu}) \dlto^{c_*}\\
H_*(Y_{\xi (\mu)}) \rrto_{{f_{\mu}}_*} && H_*(X_{\mu})
\enddiagram$$

The commutativity of the top and bottom squares is due to the naturality of the Chern--Schwartz--MacPherson class [M], the commutativity of the right and left squares is due to the above Corollary (3.9), and the commutativity of the outer big square is due to the following fact:
For the fiber square 
$$
\CD
X' @> g' >> X \\
@V f' VV @VV f V\\
Y' @>> g > Y, \endCD 
$$ 
the following diagram commutes 
(e.g., see [Er, Proposition 3.5], [FM, Axiom ($A_{23}$)]):
$$
\CD
F(Y') @> {f'}^* >>  F(X') \\
@V {g_*}VV @VV {g'_*} V\\
 F(Y) @>> f^* >  F(X). \endCD
$$
Thus it remains to see only the commutativity of the inner small square, i.e., for any $y \in H_*(Y_{\xi (\la)})$
$$\pi_{\la \mu}^{**} \Bigl ({f_{\la}}_*(y) \Bigr ) = {f_{\mu}}_* \left (\rho_{\xi (\la) \xi (\mu)}^{**} (y) \right ),$$
which is more precisely
$$\ga^{\op {Br}}(\jeden_{\pi_{\la \mu}}) \bullet_{\bH}{f_{\la}}_{\bigstar}(y) = {f_{\mu}}_{\bigstar} \Bigl (\ga^{\op {Br}} (\jeden_{\rho_{\xi (\la) \xi (\mu)}}) \bullet_{\bH} y \Bigr ).$$
Since $\jeden_{\rho_{\xi (\la) \xi (\mu)}}= f_{\la}^{\bigstar}\jeden_{\pi_{\la \mu}}$ and the Grothendieck transformation $\ga^{\op {Br}}: \bF \to \bH$ is compatible with the pullback operation, the above equality becomes
$$\ga^{\op {Br}}(\jeden_{\pi_{\la \mu}}) \bullet_{\bH}{f_{\la}}_{\bigstar}(y) = {f_{\mu}}_{\bigstar} \Bigl (f_{\la}^{\bigstar}\ga^{\op {Br}} (\jeden_{\pi_{\la \mu}}) \bullet_{\bH} y \Bigr ).$$
And it turns out that this equality is nothing but {\it the projection formula} of the Bivariant Theory [FM, \S 2.2, ($\text {A}_{\text {123}}$)] for the following diagram and for the bivariant homology theory $\bH$:
$$\CD
Y_{\xi (\mu)}@> {f_{\mu}} >> X_{\mu}\\
@V {\rho_{\xi(\la) \xi(\mu)}}VV @VV {\pi_{\la \mu}}V\\
Y_{\xi (\la)}@>> {f_{\la}} > X_{\la} @ >>{}>pt. \endCD
$$
Thus we obtain the following theorem.

\proclaim {Theorem (3.11)} Let $\{f_{\la}:  Y_{\xi(\la)} \to X_{\la} \}$ be a fiber-square pro-morphism between two pro-algebraic varieties with structure morphisms being Euler morphisms. Then we have the following commutative diagram:
$$\CD
F^{\op {pro}}(Y_{\infty}) @> {c_*^{\op {pro}}} >> H_*^{\op {pro}}(Y_{\infty})\\
@V {{f_{\infty}}_*}VV @VV {{f_{\infty}}_*}V\\
F^{\op {pro}}(X_{\infty}) @>> {c_*^{\op {pro}}} > H_*^{\op {pro}}(X_{\infty}). \endCD
$$
\qed
\endproclaim

Following the above construction, similarly we can get a proalgebraic version of the Riemann--Roch theorem,
i.e., the Baum--Fulton--MacPherson's Riemann--Roch $\tau_* : \bold K_0 \to {H_*}_{\bQ}$ constructed in [BFM]. 
\proclaim {Theorem (3.12)} Let $\{f_{\la}:  Y_{\xi(\la)} \to X_{\la} \}$ be a fiber-square pro-morphism between two pro-algebraic varieties with structure morphisms being proper local complete intersection morphisms. Then we have the following commutative diagram:
$$\CD
\bold K_0^{\op {pro}}(Y_{\infty}) @> {\tau_*^{\op {pro}}} >> H_*^{\op {pro}}(Y_{\infty})\\
@V {{f_{\infty}}_*}VV @VV {{f_{\infty}}_*}V\\
\bold K_0^{\op {pro}}(X_{\infty}) @>> {\tau_*^{\op {pro}}} > H_*^{\op {pro}}(X_{\infty}). \endCD
$$
\qed
\endproclaim

These results lead us to much more general theorems. First we introduce the following notion.
For a morphism $f:X \to Y$ and a bivariant class $b \in \bB(X \overset f \to\longrightarrow Y)$,
the pair $(f;b)$ is called a {\it bivariant-class-equipped morphism} and we just express $(f;b): X \to Y$. Let $\bigl \{(\pi_{\la \mu}; b_{\la \mu}): X_{\mu} \to X_{\la} \bigr \}$ be a system of bivariant-class-equipped morphisms. 
If a system $\bigl \{b_{\la \mu} \bigr \}$ of bivariant classes satisfies that
$$b_{\mu \nu} \bullet b_{\la \mu} = b_{\la \nu} \quad (\la < \mu < \nu) $$
then we call the system {\it a projective system of bivariant classes}, abusing words.
If $\bigl \{\pi_{\la \mu}: X_{\mu} \to X_{\la} \bigr \}$ and $\bigl \{b_{\la \mu} \bigr \}$ are projective systems, then the system $\bigl \{(\pi_{\la \mu}; b_{\la \mu}): X_{\mu} \to X_{\la} \bigr \}$ shall be called {\it a projective system of bivariant-class-equipped morphisms}.

For a bivariant theroy $\bB$ on the category $\Cal C$ and for a projective system $\bigl \{(\pi_{\la \mu}; b_{\la \mu}): X_{\mu} \to X_{\la} \bigr \}$ of bivariant-class-equipped morphisms, the inductive limit
$$\varinjlim_{\La} \Bigl \{B_*(X_{\la}),  b_{\la \mu} \bullet: B_*(X_{\la}) \to B_*(X_{\mu})\Bigr \}$$
shall be denoted by
$$B_*^{\op {pro}}\Bigl (X_{\infty}; \{b_{\la \mu} \} \Bigr )$$
emphasizing the projective system $\{b_{\la \mu} \}$ of bivariant classes, because the above inductive limit surely depends on the choice of it. For examples, in the above theorems we have that 
$$F^{\op {pro}}(X_{\infty}) = F_*^{\op {pro}}\Bigl (X_{\infty}; \bigl \{\jeden_{\pi_{\la \mu}} \bigr \} \Bigr )$$
$$\bold K_0^{\op {pro}}(X_{\infty}) = K_*^{\op {pro}}\Bigl (X_{\infty}; \bigl \{[\pi_{\la \mu}] \bigr \} \Bigr ).$$

Theorem (3.5)(i)  is generalized to the following theorem:

\proclaim {Theorem (3.13)}  Let $\Bigl \{(\pi_{n(n+1)}, \alp_{n(n+1)}): X_{n+1} \to X_n \Bigr \} $ be a projective system of bivariant-class-equipped morphisms of topologically connected algebraic varieties with 
\noindent \newline
$\alp_{n(n+1)} \in \bF(X_{n+1} \to X_n)$. And assume that the (of course constant) Euler--Poincar\'e characteristic $\chi \bigl (\pi_{n(n+1)}^{-1}(y); \alp_{n(n+1)} \bigr )$ of $\alp_{n(n+1)}$ restricted to each fiber $\pi_{n(n+1)}^{-1}(y)$ is non-zero and it shall be denoted by $\chi_f \bigl (\alp_{n(n+1)} \bigr )$. And we set $\chi_f(\alp_{01}):=1$.
Then the canonical Euler--Poincar\'e (pro)characteristic homomorphism
$$\chi^{\op {pro}}: F^{\op {pro}}\Bigl (X_{\infty}; \bigl \{\alp_{n(n+1)} \bigr \} \Bigr ) \to \bQ$$
is described by
$$\chi^{\op {pro}}\left (\sum_n [\alp_n] \right ) = \sum_n \frac {\chi(\alp_n)}{\chi_f(\alp_{01}) \cdot \chi_f(\alp_{12}) \cdots \chi_f(\alp_{(n-1)n})}. $$
\endproclaim

\demo {Proof} Let $(f, \alp) : X \to Y$ be a bivariant-class-equipped morphism of topologically connected algebraic varieties with $\alp \in \bF(X \overset f \to \longrightarrow Y)$. It follows from Proposition (3.10) that the Euler--Poincar\'e characteristic $\chi \bigl (f^{-1}(y); \alp \bigr )$ of $\alp$ restricted to each fiber $f^{-1}(y)$ is constant (and non-zero by assumption). So, if it is denoted by $\chi_f(\alp)$, then $f_*\alp = \chi_f(\alp) \cdot \jeden_Y$. Then to prove the theorem it suffices to see that we have the following commuttive diagram:
$$\CD
F(Y)@> {\alp \bullet } >> F(X) \\
@V {\chi}VV @VV {\chi}V \\
\bZ@>> {\times \chi_f(\alp)} > \bZ. \endCD
$$
To see this, we need the {\it projection formula} that for a morphism $f: X \to Y$ and constructible functions $\alp \in F(X)$ and $\be \in F(Y)$
$$f_*(\alp \cdot f^*\be) = (f_* \alp) \cdot \be.$$
Then, using this projection formula we have
$$
\align
\chi (\alp \bullet \be) & = \chi (\alp \cdot f^*\be) \\
& = \chi (f_*\alp \cdot \be) \\
& = \chi \left  ( \left (\chi_f(\alp) \cdot \jeden_Y \right ) \cdot \be \right ) \\
& = \chi_f(\alp) \cdot \chi (\be) \endalign
$$
Thus we get the above commutative diagram.
\qed
\enddemo

Theorem (3.10), Theorem (3.11) and Theorem (3.12) are generalized to the following theorem:

\proclaim {Theorem (3.14)} (i) Let $\ga: \bB \to \bB'$ be a Grothendieck transformation between two bivariant theories $\bB, \bB': \Cal C \to \Cal C'$ and let $\bigl \{(\pi_{\la \mu}; b_{\la \mu}): X_{\mu} \to X_{\la}) \bigr \}$ be a projective system of bivariant-class-equipped morphisms. 
Then we get the following pro-version of the natural transformation $\ga_*:B_* \to B'_*$:
$$\ga_*^{\op {pro}}: B_*^{\op {pro}}\Bigl (X_{\infty}; \{b_{\la \mu} \} \Bigr ) \to {B'_*}^{\op {pro}}\Bigl (X_{\infty}; \{\ga (b_{\la \mu}) \} \Bigr ).$$

(ii)  Let $\{f_{\la}:  Y_{\xi(\la)} \to X_{\la} \}$ be a fiber-square pro-morphism between two projective systems of bivariant-class-equipped morphisms such that $b_{\xi(\la) \xi(\mu)}= f_{\la}^{\bigstar}b_{\la \mu}$.
Then we have the following commutative diagram:
$$\CD
B_*^{\op {pro}}(Y_{\infty}) @> {\ga_*^{\op {pro}}} >> {B'}_*^{\op {pro}}(Y_{\infty})\\
@V {{f_{\infty}}_*}VV @VV {{f_{\infty}}_*}V\\
B_*^{\op {pro}}(X_{\infty}) @>> {\ga_*^{\op {pro}}} > {B'}_*^{\op {pro}}(X_{\infty}). \endCD
$$
\qed
\endproclaim

We hope to be able to do further investigations on characteristic classes (in particular, characteristic classes having a bivariant version) of proalgebraic varieties and some applications of them.

\head \S 4 $\chi$-stable proconstructible functions \endhead

In the previous section we have dealt with the case when the Euler--Poincar\'e characteristic of the fibers of each structure (or bonding) morphism, or the Euler--Poincar\'e characteristic of a constructible function restricted to the fibers of each structure morphism is non-zero constant.
In this section we address ourselves to more general cases when this condition does not necessarily hold.

Let $\bigl \{p_{\la \mu} \bigr \}$ be a system of non-zero integers indexed by the directed set $\La$. If the following holds, then $\bigl \{p_{\la \mu} \bigr \}$ shall be called a projective system:
$$p_{\la \la} = 1 \qquad \text {and} \qquad p_{\la \mu} \cdot p_{\mu \nu}  = p_{\la \nu} \quad (\la < \mu < \nu) .$$
For each $\la \in \La$ we define the following subgroup of $F(X_{\la})$:
$$F^{\op {st}}_{\{p_{\la \mu} \}}(X_{\la}) : = \Bigl \{ \alp_{\la} \in F(X_{\la}) \quad | \quad 
\chi \bigl (\pi_{\la \mu}^*\alp_{\la} \bigr ) 
= p_{\la \mu}\cdot \chi(\alp_{\la}) \quad \text {for any} \quad \mu > \la \Bigr \}.$$
For each $\la \in \La$, an element of $F^{\op {st}}_{\{p_{\la \mu} \}}(X_{\la})$ is called 
{\it a $\chi$-stable constructible function with respect to the projective system $\bigl \{p_{\la \mu} \bigr \}$ of non-zero integers }. Then it is easy to see that for each structure morphism $\pi_{\la \mu} : X_{\mu} \to X_{\la}$ the pullback homomorphism $\pi_{\la \mu}^*: F(X_{\la}) \to F(X_{\mu})$ preserves  $\chi$-stable constructible functions with respect to the projective system $\bigl \{p_{\la \mu} \bigr \}$ of non-zero integers, namely it induces the homomorphism (using the same symbol):
$$\pi_{\la \mu}^*: F^{\op {st}}_{\{p_{\la \mu} \}}(X_{\la}) \to F^{\op {st}}_{\{p_{\la \mu} \}}(X_{\mu})$$
which implies that we get the inductive system
$$\Bigl \{F^{\op {st}}_{\{p_{\la \mu} \}}(X_{\la}) , \pi_{\la \mu}^*: F^{\op {st}}_{\{p_{\la \mu} \}}(X_{\la}) \to F^{\op {st}}_{\{p_{\la \mu} \}}(X_{\mu}) \quad (\la < \mu) \Bigr \}.$$ 
Then for a proalgebraic variety $\displaystyle X_{\infty} = \varprojlim _{\la \in \La} X_{\la}$  we consider the inductive limit of the above inductive system and it shall be denoted by
$$F^{\op {st.pro}}_{\{p_{\la \mu} \}}(X_{\infty})$$
and an element of this group shall be called a {\it $\chi$-stable proconstructible function on the  proalgebraic variety $\displaystyle X_{\infty}$ with respect to the projective system $\bigl \{p_{\la \mu} \bigr \}$ of non-zero integers}.
We can see that this subgroup can be also directly defined as follows:
$$ \Bigl \{ [\alp_{\la}] \in F^{\op {pro}}(X_{\infty}) \quad | \quad \chi (\pi_{\la \mu}^*\alp_{\la}) 
= p_{\la \mu}\cdot \chi(\alp_{\la}) \quad (\la < \mu) \Bigr \}.$$

For each structure morphism $\pi_{\la \mu} : X_{\mu} \to X_{\la}$ we get the following commutative diagram

$$\CD
F^{\op {st}}_{\{p_{\la \mu} \}}(X_{\la})@> {\pi_{\la \mu}} >>F^{\op {st}}_{\{p_{\la \mu} \}}(X_{\mu}) \\
@V {\chi}VV @VV {\chi}V \\
\bZ@>> {\times p_{\la \mu}} > \bZ. \endCD
$$

Thus we can get the following theorem:

\proclaim {Theorem (4.1)} For a proalgebraic variety $\displaystyle X_{\infty}= \varprojlim _{\la \in \La} X_{\la}$ and a projective system $\bigl \{p_{\la \mu} \bigr \}$ of non-zero integers, we get the proalgebraic Euler--Poincar\'e characteristic homomorphism
$$\chi^{\op {st.pro}}_{\{p_{\la \mu} \}}: F^{\op {st.pro}}_{\{p_{\la \mu} \}}(X_{\infty}) \to \varprojlim _{\la \in \La} \Bigl \{ \times p_{\la \mu} : \bZ \to \bZ \Bigr \}.$$
\endproclaim

In particular, when it comes to the case when $\La = \bN$, we get the following theorem
\proclaim {Theorem (4.2)}  For a proalgebraic variety $\displaystyle X_{\infty} = \varprojlim _{n \in \bN} X_n$ and a projective system $\{p_{n m} \}$ of non-zero integers, we have the following
canonical proalgebraic Euler--Poincar\'e characteristic homomorphism
$$\chi^{\op {st.pro}}_{\{p_{n m} \}}: F^{\op {st.pro}}_{\{p_{n m} \}}(X_{\infty}) \to \bQ$$
which is defined by
$$\chi^{\op {st.pro}}_{\{p_{n m} \}}\Bigl (\sum_n [\alp_n] \Bigr ) 
:= \sum_n \frac {\chi(\alp_n)}{p_{0 1}\cdot p_{1 2}\cdot p_{2 3}\cdots p_{(n-1) n}}.$$
Here we set $p_{0 1}:= 1$.
\endproclaim

In the previous section all proconstructible functions are $\chi$-stable with respect to some projective systems of non-zero integers.

\head \S 5 motivic measures \endhead

The theorems or observations obtained concerning Chern-Schwartz-MacPhrson classes or Euler--Poincar\'e characteristics of proalgebraic varieties in the previous section led the author to realize some deep connections to motivic measures and motivic integrations, in which the Grothendieck ring of complex proalgebraic varieties is a crucial ingredient.

First we recall the usual Grothendieck ring of algebraic varieties.
Let $\Cal V_{\bC}$ denote the category of complex algebraic varieties. 
Then the Grothendieck ring $K_0(\Cal V_{\bC})$ of complex algebraic 
varieties is the free abelian group generated by the isomorphism classes of varieties
 modulo the subgroup generated by elements of the form $[V] - [V'] - [V \setminus V']$ 
for a closed subset $V' \subset V$ with the ring structure $[V]\cdot [W] := [V \times W]$.
 There are distinguished elements in $K_0(\Cal V_{\bC})$: $\jeden$ is the class [p] of a point $p$
 and $\bL$ is the Tate class $[\bC]$ of the affine line $\bC$. From this definition, 
we can see that any constructible set of a variety determines an element 
in the Grothendieck ring $K_0(\Cal V_{\bC})$. Provisionally the element $[V]$
 in the Grothendieck ring $K_0(\Cal V_{\bC})$ is called the {\it Grothendieck ``motivic" class of $V$} 
and let us denote it by $\Ga (V)$. Hence we get the following homomorphism, 
called the {\it Grothendieck ``motivic" class homomorphism}: for any variety $X$
$$\Ga : F(X) \to K_0(\Cal V_{\bC}),$$
which is defined by
$$\Ga (\alp) = \sum _{n \in \bZ} n \left [\alp^{-1}(n) \right ].$$
Or $\Ga \left (\sum a_V \jeden_V \right ) := \sum a_V [V]$
where $V$ is a constructible set in $X$ and $a_V \in \bZ$.

This Grothendieck ``motivic" class homomorphism is {\it tautological} and its more ``geometric" one is the Euler--Poincar\'e characteristic homomorphsim $\chi: F(X) \to \bZ$. In the previous section we have generalized $\chi: F(X) \to \bZ$ to the category of proalgebraic varieties.
Thus our problem is to generalize the Grothendieck ``motivic" class homomorphism $\Ga : F(X) \to K_0(\Cal V_{\bC})$ to the category of proalgebraic varieties. 

\proclaim {Problem (5.1)} Let $\displaystyle X_{\infty}= \varprojlim_{\la \in \La} X_{\la}$ be a proalgebraic variety. Then describe a proalgebraic Grothendieck ``motivic" (pro)class homomorphism
$$\Ga^{\op {pro}} : F(X_{\infty}) \to \sharp$$
with some reasonable algebraic object $\sharp$.
\endproclaim

It seems that one cannot expect a general formula.
However, we can get a Grothendieck ``motivic" class homomorphism version of Theorem (3.5), if we consider the projective limit of a projective sytem of Zariski (not in the usual topology) locally trivial fiber bundles. Indeed, the Grothendieck ``motivic" (pro)class homomorphism version of (3.5.2) is the following:

$$
\CD 
F(X_n) @ > {\pi_{(n-1)n}^*} >> F(X_{n+1}) \\
@V {\Ga} VV  @ VV {\Ga} V \\
K_0(\Cal V_{\bC}) @>> {\times [F_n]} > K_0(\Cal V_{\bC}). \endCD \tag 5.2
$$

However, when we consider a localization of the Grothendieck ring $K_0(\Cal V_{\bC})$, we need to be a bit careful. Namely, the Grothendieck ring $K_0(\Cal V_{\bC})$ is not a domain unlike the ring $\bZ$ of integers. Indeed, here is a very recent result due to B. Poonen [Po, Theorem 1]:
\proclaim {Theorem (5.3)}(B. Poonen) Suppose that $k$ is a field of characteristic zero. Then the Grothendieck ring $K_0(\Cal V_k)$ of $k$-varieties is not a domain.
\endproclaim

Recall that for any ring $R$ and $S$ the multplicatively closed set of all non-zeo divisors of $R$, the quotient ring $R_S$ of $R$ with denominator set being $S$ is called the {\it total (or full) quotient ring of $R$} and denoted by $Q(R)$ (eg. see [Kunz]). So, let $Q(K_0(\Cal V_{\bC}))$ denote the total quotient ring of the Grothendieck ring $K_0(\Cal V_{\bC})$. Then, very simply thinking, the ``motivic" proclass version of the canonical Euler--Poincar\'e (pro)characteristic homomorphism $\chi^{\op {pro}}: F^{\op {pro}}(X_{\infty}) \to \bQ$ given in Theorem (3.5)
could be the homomorphism 
$$\Ga^{\op {pro}}: F^{\op {pro}}(X_{\infty}) \to Q(K_0(\Cal V_{\bC}))$$
defined by
$$\Ga^{\op {pro}}\left (\sum_n [\alp_n] \right ) = \sum_n \frac {\Ga (\alp_n)}{[F_0] [F_1] [F_2] \cdots [F_{n-1}]}. $$
Here $F_0$ is defined to be a point, i.e.,$[F_0]:= \jeden$, and also we assume that each $[F_j]$ is not a zero-divisor.

One serious problem of this definition is that it is in general hard to check whether the class $[F_j]$ is a non-zero divisor or not; indeed, we do not know whether even the class $[\bP^j]$ of the projective space is a non-zero divisor (which Willem Veys pointed out to the author). So, in order to avoid this problem, we consider the following simple localization: let $\Cal F$ be the multiplicative set consisting of all the finite products of $[F_j]^{m_j}$, i.e, 
$$\Cal F := \Bigl \{ [F_{j_1}]^{m_1}[F_{j_2}]^{m_2}\cdots [F_{j_s}]^{m_s} | j_i \in \bN, m_i \in \bN \Bigr \}.$$

Then the quotient ring $\Cal F^{-1}\left (K_0(\Cal V_{\bC}) \right )$ shall be denoted by $K_0(\Cal V_{\bC})_\Cal F$. With this definition we get the following theorem:

\proclaim {Theorem (5.4)} Let $\displaystyle X_{\infty}= \varprojlim_n X_n$ be a proalgebraic variety such that each structure morphism $\pi_{n(n+1)}: X_{n+1} \to X_n$ satisfies the condition that $[X_{n+1}] = [X_n][F_n]$ for some variety $F_n$ for each $n$, for example, $\pi_{n, n+1}: X_{n+1} \to X_n$ is a Zariski locally trivial fiber bundle with fiber variety being $F_n$.

(i) The canonical Grothendieck ``motivic" proclass homomorphism,
$$\Ga^{\op {pro}}: F^{\op {pro}}(X_{\infty}) \to K_0(\Cal V_{\bC})_\Cal F$$
is described by
$$\Ga^{\op {pro}}\left (\sum_n [\alp_n] \right ) = \sum_n \frac {\Ga (\alp_n)}{[F_0] [F_1] [F_2] \cdots [F_{n-1}]}. $$
Here $F_0$ is defined to be a point, i.e.,$[F_0]:= \jeden$. 

(ii) In particular, if all the fibers are the same, say $F_n = W$ for any $n$, then the canonical Grothendieck ``motivic" (pro)class homomorphism
$$\Ga^{\op {pro}}: F^{\op {pro}}(X_{\infty}) \to K_0(\Cal V_{\bC})_\Cal F$$
is described by
$$\Ga^{\op {pro}}\left (\sum_n [\alp_n] \right ) = \sum_n \frac {\Ga (\alp_n)}{[W]^{n-1}}. $$
In this special case the quotient ring $K_0(\Cal V_{\bC})_\Cal F$ shall be simply denoted by $K_0(\Cal V_{\bC})_{[W]}$.
And another non-canonical one $\Ga_w^{\op {pro}}$ is 
$$\Ga_w^{\op {pro}}\left (\sum_n [\alp_n] \right ) = \sum_n \frac {\Ga (\alp_n)}{[W]^{n-1+w }}. $$
 \qed
\endproclaim

Let $W$ be a variety and let $X \widetilde {\times} W^k$ be a fiber bundle inductively defined as a fiber bundle over $X \widetilde {\times} W^{k-1}$ with fiber being $W$. Then the projective limit $\displaystyle \varprojlim_k X \widetilde {\times} W^k$ shall be denoted by $X \widetilde {\times} W^{\bN}$. Then, for the characteristic function $\jeden_{X \widetilde {\times} W^{\bN}}=  [\jeden_{X \widetilde {\times} W^{k-1}}]$ for any $k$, we have 
$$\Ga^{\op {pro}}\left (\jeden_{X \widetilde {\times} W^{\bN}}\right ) = \frac{[X]\cdot [W]^{k-1}}{[W]^{k-1}} = [X],$$
which is called {\it the Grothendieck ``motivic" {\it proclass} of the proalgebraic variety $X \widetilde {\times} W^{\bN}$} and denoted by $[X \widetilde {\times} W^{\bN}]^{\op {pro}}$ and sometimes called a ``fiber bundle" over an algebraic variety $X$ with fiber being the proalgebraic variety $W^{\bN}$, abusing words. 

As to the Grothendieck class ``motivic" homomorphism version of Theorem (3.5)(iii), we do not know whether it is true or not. Thus we want to ask the following question:

\proclaim {Question (5.5)} Is the canonical Grothendieck ``motivic" (pro)class homomorphism
$$\Ga^{\op {pro}}: F^{\op {pro}}(X_{\infty}) \to K_0(\Cal V_{\bC})_\Cal F$$
is surjective if each $F_n = \bP^n$ the projective space of dimension $n$ ?
\endproclaim

Since the Grothendieck ``motivic" (pro)classes of proconstructible functions on the above proalgebraic variety $X \widetilde {\times} W^{\bN}$ are in the localization $K_0(\Cal V_{\bC})_{[W]}$ of the Grothendieck ring $K_0(\Cal V_{\bC})$, our very na\" \i ve question is:
\proclaim {Question (5.6)} What is the Grothendieck ring of the category $\Cal V_{\bC}^{\op {pro}}$ of proalgebraic varieties ? What is $K_0(\Cal V_{\bC}^{\op {pro}})$ ?
\endproclaim
At the moment we do not have an answer for this question. However, we can give a plausible answer for this in a special category of proalgebraic varieties, which are the above ``fiber bundles" over complex algebraic varieties with fiber being the proalgebraic variety $W^{\bN}$. Infinite ``fiberings with fiber $W$" give rise to the following  inductive system
$$ \cdots \overset {\times [W]} \to \longrightarrow K_0(\Cal V_{\bC}) \overset {\times [W]} \to \longrightarrow  K_0(\Cal V_{\bC}) \overset {\times [W]} \to \longrightarrow  K_0(\Cal V_{\bC}) \overset {\times [W]} \to \longrightarrow \cdots $$
the colimit of which would be the localization $K_0(\Cal V_{\bC})_{[W]}$ of the Grothendieck ring $K_0(\Cal V_{\bC})$.

Now we discuss motivic measures of the arc space $\Cal L(X)$ of an algebraic variety $X$.
The arc space $\Cal L(X)$ of an algebraic variety $X$ is
 defined to be the projective limit of truncated arc varieties $\Cal L_n(X)$ and projection $\pi_{n-1, n}: \Cal L_n(X) \to \Cal L_{n-1}(X)$. Thus the arc space $\Cal L(X)$ is a proalgebraic
 variety. Let $\pi_n: \Cal L(X) \to \Cal L_n(X) $ be the canonical projection or
the $n$-th truncated morphism. Then a subset $A$ of the arc space $\Cal L(X)$ is called 
{\it a cylinder set} [Cr] or simply {\it a constructible set} [DL 1, DL 2] if $A = \pi_n^{-1}(C_n)$ for a constructible set $C_n$ in the $n$-th arc space $\Cal L_n(X)$ for some integer $n \geq 0$. 
In this paper, to avoid some possible confusion, we take the name ``cylinder set". 

A $\bZ$-valued function $\alp:\Cal L(X) \to \bZ$ on the arc space $\Cal L(X)$ is called a {\it cylinder function} if $\alp$ is a linear combination of characteristic functions on cylinder sets. 
It turns out that the abelian group of cyclinder functions on the arc space $\Cal L(X)$, denoted by $F^{\op {cyl}}(\Cal L(X))$, is isomorphic to the abelian group $F^{\op {pro}}(\Cal L(X))$, as shown below.

Our previous definition of a proconstructible function in $F^{\op {pro}}(X_{\infty}) = \varinjlim _{\la \in \La} F(X_{\la})$ as a function on $X_{\infty}$ means that we define the ``functionization" homomorphism 
$$\Psi : \varinjlim _{\la \in \La} F(X_{\la}) \to Fun(X_{\infty}, \bZ) \quad \text {defined by} \quad
\Psi \left ([\alp_{\mu}] \right ) \left ((x_{\la}) \right ) := \alp_{\mu}(x_{\mu}).$$
Here $Fun(X_{\infty}, \bZ)$ is the abelian group consisting of functions or mappings from $X_{\infty}$ to  $\bZ$. Of course, the target $\bZ$ can be replaced by any other ring containing the integers $\bZ$.

One can describe the above in a much fancier way as follows.  Let $\pi_{\la}: X_{\infty} \to X_{\la}$ denote the canonical projection induced from the projection $\prod_{\la} X_{\la} \to X_{\la}$. Consider the following commutaive diagram for each $\la \in \La$:

$$\xymatrix{
F(X_\lambda) \ar[dd]_{\pi^*_{\la \mu}} \ar[dr] ^{\pi^*_\la} \\
& Fun(X_\infty,\Bbb Z) \\
F(X_\mu) \ar[ur]_{\pi^*_\mu}
}$$

Then it follows as a standard fact in the theory of inductive limits that the ``functionization" homomorphism $\Psi : \varinjlim _{\la \in \La} F(X_{\la}) \to Fun(X_{\infty}, \bZ)$ is the unique homomorphism such that the following diagram commutes:

$$\xymatrix{
& F(X_{\la}) \ar [dl]_{ \rho^{\la}} \ar [dr]^{\pi_{\la}^*} \\
{\varinjlim _{\la \in \La} F(X_{\la}) } \ar [rr] _{\Psi}& &  {Fun(X_{\infty}, \bZ)}.}
$$

To avoid some possible confusion, the image $\Psi \bigl ([\alp_{\la}] \bigr ) = \pi_{\la}^*\alp_{\la}$ shall be denoted by $[\alp_{\la}]_{\infty}$.   For a constructible set $W_{\la} \in X_{\la}$, $[\jeden_{W_{\la}}]_{\infty}$ shall be also called a {\it procharacteristic function}.  Then $\op {Image}\Psi$ is generated by all the procharacteristic functions $[\jeden_{W_{\la}}]_{\infty}$ with $W_{\la} \subset X_{\la}$ constructible sets. By the definition we have
$$[\jeden_{W_{\la}}]_{\infty} = \jeden_{\pi_{\la}^{-1}(W_{\la})}.$$
Thus, mimicking the term of cylinder set [Cr], the set $\pi_{\la}^{-1}(W_{\la}) \subset X_{\infty}$ for a constructible set $W_{\la} \in X_{\la}$ is also called {\it a $\la$-cylinder set} or {\it a cylinder set of level $\la$}. Thus $[\jeden_{W_{\la}}]_{\infty}$ is a characteristic function on the cylinder set, so a cylinder function. It is easy to see the following

\proclaim {Proposition (5.7)} For $\la < \mu$ and two constructible sets $W_{\la} \in X_{\la}$ and
$W_{\mu} \in X_{\mu}$ we have

$$\pi_{\la}^{-1}(W_{\la}) = \pi_{\mu}^{-1}(W_{\mu}) 
\Longleftrightarrow W_{\mu}= \pi_{\la \mu}^{-1} \bigl (W_{\la}\bigr )
\Longleftrightarrow \jeden_{W_{\mu}} =  \pi_{\la \mu}^* \Bigl (\jeden_{W_{\la}} \Bigr ).$$
\endproclaim

Thus we can see that the notions of procharacteristic function and of cylinder set are {\it  equivalent}. 

Thus we define

\definition {Definition (5.8)} For a proalgebraic variety $\displaystyle X_{\infty} = \varprojlim _{\la \in \La} X_{\la}$
$$F^{\op {cyl}}(X_{\infty}) := \op {Image}\left (\Psi : \varinjlim _{\la \in \La} F(X_{\la}) \to Fun (X_{\infty}, \bZ) \right ) = \bigcup_{\mu} \pi_{\mu}^* \bigl  (F(X_{\mu}) \bigr ).$$
And an element of $F^{\op {cyl}}(X_{\infty})$ is called a cylinder function on $X_{\infty}$.
\enddefinition

\remark {Remark (5.9)} As long as the support is concerned,

$$\op {supp}[\jeden_{W_{\la}}] = \op {supp}[\jeden_{W_{\la}}]_{\infty} = \pi_{\la}^{-1}(W_{\la}),$$
whether the characteristic function $\jeden_{W_{\la}}\in F(X_{\la})$ is considered as the proconstructible function $[\jeden_{W_{\la}}] \in F^{\op {pro}}(X_{\infty})$ or considered as the cylinder function $[\jeden_{W_{\la}}]_{\infty} \in F^{\op {cyl}}(X_{\infty})$. Note that our results obtained so far of course hold with nothing changed at all even if $F^{\op {pro}}(X_{\infty})$ is replaced by $F^{\op {cyl}}(X_{\infty})$.
\endremark

In an earlier version of the present paper, we considered only $F^{\op {pro}}(X_{\infty})$. However, in a discussion with J\"org Sch\"urmann, he suggested to also consider the image of the above ``functionization" homomorphism, pointing out that in general the inductive limit would have more informations than the image, i.e., that $\Psi$ may not be necessarily injective.

In some cases it is injective. For example,  if all the homomorphism $\pi_{\la}^*: F(X_{\la}) \to Fun(X_{\infty}, \bZ)$ are injective, which is in turn equivalent to the condition that all the projections $\pi_{\la}:X_{\infty} \to X_{\la}$ is surjective, then the ``functionization" homomorphism $\Psi : \varinjlim _{\la \in \La} F(X_{\la}) \to Fun(X_{\infty}, \bZ)$ is also injective since the inductive limit is an exact functor. So, we get the following lemma:

\proclaim {Proposition (5.10)} If all the structure morphisms $\pi_{\mu \la}: X_{\mu} \to X_{\la}$ (for $\mu < \la$) are surjective, then the ``functionization" homomorphism $\Psi : \varinjlim _{\la \in \La} F(X_{\la}) \to Fun(X_{\infty}, \bZ)$ is injective. 
\endproclaim
 
In the case of the arc space $\Cal L(X)$, since each structure morphism $\pi_{n, n+1}: \Cal L_{n+1}(X) \to \Cal L_n (X)$ is always surjective, we get the following 

\proclaim {Corollary (5.11)} For the arc space $\Cal L(X)$  we have the canonical isomorphism
$$F^{\op {pro}}\bigl (\Cal L(X) \bigr )\cong F^{\op {cyl}}\bigl (\Cal L(X) \bigr ).$$
\endproclaim

Suppose that $\Psi ([\alp_{\mu}]) = 0$, which means that $\Psi ([\alp_{\mu}]) ((x_{\la}) = \alp_{\mu}(x_{\mu}) = 0$ for any $(x_{\la}) \in X_{\infty}$. Hence we have
$$\alp_{\mu}\bigl (\pi_{\mu}(X_{\infty}) \bigr ) = 0.$$
At the moment we do not know whether we can conclude $[\alp_{\mu}] = 0$ from this condition.
There is a very simple example such that $\alp_{\mu}\bigl (\pi_{\mu}(X_{\infty}) \bigr ) = 0, \pi_{\mu}(X_{\infty}) \not = X_{\mu}$ and $\alp_{\mu} \not = 0$, but $[\alp_{\mu}] = 0:$
Let $X_1 = \{a, b \}$ be a space of two different points, and let $X_n = \{a \}$ for any $n > 1$.
Let $\pi_{1, 2}:X_2 \to X_1$ is the injection map sending $a$ to $a$ and the other structure
morphism $\pi_{n(n+1)}:X_{n+1} \to X_n$ is the identity for $n>1$. Then the projective limit
$X_{\infty} = \{ (a) \}$ consits of one point $(a, a, a, \cdots)$. Let $\alp_1 = p\cdot \jeden_b \in F(X_1)$. Then we have $\alp_1\left (\pi_1(X_{\infty}) \right ) = 0, \pi_1(X_{\infty}) \not = X_1$ and $\alp_1 \not = 0$, but $[\alp_1] = 0.$ We suspect that in general the ``functionization" homomorphism $\Psi$ is not necessarily injective, but we have not been able to find such an example yet.
 
\proclaim {Corollary (5.12)} When $X$ is a nonsingular variety of dimension $d$, we have the following canonical Grothendieck ``motivic" (pro)class homomorphism
$$\Ga^{\op {pro}}: F^{\op {cyl}}(\Cal L(X)) \to K_0(\Cal V_{\bC})_{[\bL^d]}$$
is described by
$$\Ga^{\op {pro}}\left (\sum_n [\alp_n]_{\infty} \right ) = \sum_n \frac {\Ga (\alp_n)}{[\bL]^{nd}} $$
and another non-canonical one $\Ga_w^{\op {pro}}$ is
$$\Ga_w^{\op {pro}}\left (\sum_n [\alp_n]_{\infty} \right ) = \sum_n \frac {\Ga (\alp_n)}{[\bL]^{(n+w)d}} $$
In particular, we get that $[\Cal L(X)]^{\op {pro}}= [\jeden_{\Cal L(X)}]^{\op {pro}} = [X].$
\qed
\endproclaim
Note that in the case of arc space $\Cal L(X)$, since $\Cal L_0(X) = X$, the index set is not $\bN$ but $\{0\}\cup \bN$.  Hence the canonical one is not $\displaystyle \Ga^{\op {pro}}\left (\sum_n [\alp_n]_{\infty} \right ) = \sum_n \frac {\Ga (\alp_n)}{[\bL]^{(n-1)d}}$.

Therefore we can see that our proalgebraic Grothendieck ``motivic" (pro)class homomorphism $\Ga: F^{\op {pro}}(X_{\infty}) \to Q \left (K_0(\Cal V_{\bC}) \right )$ is a generalization of the so-called motivic measure in the case when the base variety $X$ is smooth. 

Thus the theorem and corollary imply that in the theory of motivic measures and motivic integrations the notion of cylinder set is a right one from the proalgebraic viewpoint or from the viewpoint of pro-category.

If $X$ is singular, the arc space $\Cal L(X)$ is {\it not} the projective limit of a projective system of Zariski locally trivial fiber bundles with fiber being $\bC^{\op {dim}X}$ any longer and each projection morphism $\pi_{(n-1)n}: \Cal L_n(X) \to \Cal L_{n-1}(X)$ is complicated and thus as a proalgebraic variety $\Cal L(X)$ is complicated and in general we do not know whether there exists a well-defined
homomorphism $\Ga: F^{\op {pro}}(\Cal L(X)) \to Quot(K_0(\Cal V_{\bC}))$ for some suitable quotient ring $Quot(K_0(\Cal V_{\bC}))$ , also we do not know whether there exists a well-defined proalgebraic Euler--Poincar\'e characteristic
homomorphism $\chi^{\op {pro}}: F^{\op {pro}}(\Cal L(X)) \to \bQ$.

A crucial ingredient in studing motivic measure or motivic integration is the so-called 
{\it stable set} of the arc space $\Cal L(X)$. 
\definition {Definition (5.13)} A subset $A$ of the arc space $\Cal L(X)$ is called {\it a stable set} if it is a cylinder set,
i.e., $A = \pi_n^{-1}(C_n)$ for a constructible set $C_n$ in the $n$-th arc space $\Cal L_n(X)$,  such that the restriction of each projection $\pi_{m(m+1)}|_{\pi_{m+1}(A)}: \pi_{m+1}(A) \to \pi_m(A)$ for  each $m \geq n$ is a Zariski locally fiber bundle with the fiber being $\bC^{\op {dim}X}$, in other words, the restriction of the projection $\pi_n: \Cal L(X) \to \Cal L_n(X)$ to the constructible set $C_n$ is $C_n \widetilde {\times}(\bC^{\op {dim}X})^{\bN}$. And a {\it stable function} on the arc space $\Cal L(X)$ is a $\bZ$-valued function constant along stable sets. 
\enddefinition

In fact, in a similar way as Definition (5.13) or as in \S 4, we can define a {\it ``motivic" stable proconstructible function} on any proalgebraic variety $X_{\infty}$ as follows.

Let $\bigl \{[F_{\la \mu}] \bigr \}$ be a system of Grothendieck classes $[F_{\la \mu}] \in K_0(\Cal V_{\bC})$  indexed by the directed set $\La$. As in \S 4, if the following holds, then $\bigl \{[F_{\la \mu}] \bigr \}$ shall be called a projective system:
$$[F_{\la \la}] = \jeden \qquad \text {and} \qquad [F_{\la \mu}] \cdot [F_{\mu \nu}]  = [F_{\la \nu}] \quad (\la < \mu < \nu). $$
For each $\la \in \La$ we define the following subgroup of $F(X_{\la})$:

$$F^{\op {st}}_{\{[F_{\la \mu}] \}}(X_{\la}) : = \Bigl \{ \alp_{\la} \in F(X_{\la}) \quad | \quad 
\Ga \bigl (\pi_{\la \mu}^*\alp_{\la} \bigr ) 
= [F_{\la \mu}]\cdot \Ga(\alp_{\la}) \quad \text {for any} \quad \mu > \la \Bigr \}.$$

For each $\la \in \La$, an element of $F^{\op {st}}_{\{[F_{\la \mu}] \}}(X_{\la})$ is called 
{\it a $\Ga$-stable constructible function with respect to the projective system $\bigl \{[F_{\la \mu}] \bigr \}$ of Grothendieck classes }. Then it is easy to see that for each structure morphism $\pi_{\la \mu} : X_{\mu} \to X_{\la}$ the pullback homomorphism $\pi_{\la \mu}^*: F(X_{\la}) \to F(X_{\mu})$ preserves  $\Ga$-stable constructible functions with respect to the projective system $\bigl \{[F_{\la \mu}] \bigr \}$ of Grothendieck classes, namely it induces the homomorphism (using the same symbol):

$$\pi_{\la \mu}^*: F^{\op {st}}_{\{[F_{\la \mu}]  \}}(X_{\la}) \to F^{\op {st}}_{\{[F_{\la \mu}]  \}}(X_{\mu})$$
which implies that we get the inductive system

$$\Bigl \{F^{\op {st}}_{\{[F_{\la \mu}]  \}}(X_{\la}) , \pi_{\la \mu}^*: F^{\op {st}}_{\{[F_{\la \mu}]  \}}(X_{\la}) \to F^{\op {st}}_{\{[F_{\la \mu}]  \}}(X_{\mu}) \quad (\la < \mu) \Bigr \}.$$ 
Then for a proalgebraic variety $\displaystyle X_{\infty} = \varprojlim _{\la \in \La} X_{\la}$  we consider the inductive limit of the above inductive system and it is denoted by
$$F^{\op {st.pro}}_{\{[F_{\la \mu}]  \}}(X_{\infty})$$
and an element of this group shall be called a {\it $\Ga$-stable proconstructible function on the  proalgebraic variety $\displaystyle X_{\infty}$ with respect to the projective system $\bigl \{[F_{\la \mu}]  \bigr \}$ of Grothendieck classes}.
We can see that this subgroup can be also directly defined as follows:
$$ \Bigl \{ [\alp_{\la}] \in F^{\op {pro}}(X_{\infty}) \quad | \quad \Ga (\pi_{\la \mu}^*\alp_{\la}) 
= [F_{\la \mu}]\cdot \Ga(\alp_{\la}) \quad (\la < \mu) \Bigr \}.$$

For each structure morphism $\pi_{\la \mu} : X_{\mu} \to X_{\la}$ we get the following commutative diagram

$$\CD
F^{\op {st}}_{\{[F_{\la \mu}] \}}(X_{\la})@> {\pi_{\la \mu}} >>F^{\op {st}}_{\{[F_{\la \mu}] \}}(X_{\mu}) \\
@V {\Ga}VV @VV {\Ga}V \\
K_0(\Cal V_{\bC})@>> {\times [F_{\la \mu}]} > K_0(\Cal V_{\bC}). \endCD
$$
Thus we can get the following theorem:

\proclaim {Theorem (5.14)} For a proalgebraic variety $\displaystyle X_{\infty}= \varprojlim _{\la \in \La} X_{\la}$ and a projective system $\bigl \{[F_{\la \mu}] \bigr \}$ of Grothendieck classes, we get the Grothendieck ``motivic" class homomorphism

$$\Ga^{\op {st.pro}}_{\{[F_{\la \mu}] \}}: F^{\op {st.pro}}_{\{[F_{\la \mu}] \}}(X_{\infty}) \to \varprojlim _{\la \in \La} \Bigl \{ \times [F_{\la \mu}] : K_0(\Cal V_{\bC}) \to K_0(\Cal V_{\bC}) \Bigr \}.$$
\endproclaim

In particular we get the following theorem, which is a``motivic" version of Theorem (4.1):
\proclaim {Theorem (5.15)}  For a proalgebraic variety $\displaystyle X_{\infty} = \varprojlim _{n \in \bN} X_n$ and a projective system $\{[F_{n,m}] \}$ of Grothendieck classes, we have the following
canonical Grothendieck ``motivic" class  homomorphism
$$\Ga^{\op {st.pro}}_{\{[F_{n m}] \}}: F^{\op {st.pro}}_{\{[F_{n m}] \}}(X_{\infty}) \to K_0(\Cal V_{\bC})_\Cal F$$
which is defined by
$$\Ga^{\op {st.pro}}_{\{[F_{n m}] \}}\Bigl (\sum_n [\alp_n] \Bigr ) 
:= \sum_n \frac {\Ga(\alp_n)}{[F_{0 1}] [F_{1 2}] [F_{2 3}]\cdots [F_{(n-1) n}]}.$$
Here we set $[F_{0 1}]:= \jeden$ and $\Cal F$ is the multiplicative set consisting of all the finite products of $[F_{j (j+1)}]^{m_j}$ as in Theorem (5.4).
\endproclaim

In the case of the arc space $\Cal L(X)$, the projective system $\{[F_{n m}] \}$ of Grothendeick classes is such that $[F_{n m}] = [\bC^{(m-n)d}].$ Thus Theorem (5.14) is a generalization of motivic measure when we deal with the arc space $\Cal L(X)$ of a possibly singular variety $X$.

\remark {Remark (5.16)} So far we have looked at the two components of $c_*:F(X) \to H_*(X)$. 
So it is natural to consider whether or not there is a motivic version of the whole $c_*$, say,
$$\mu c_*: F(X) \to K_0(\Cal V_{\bC})$$
such that the two distinguished ``components" of $\mu c_*$ are 
$$\chi : F(X)  \to \bZ \qquad \text {and} \qquad \Ga: F(X) \to K_0(\Cal V_{\bC}).$$
A very na\"\i ve and very simple-minded guess is the following ``construction":  Suppose that for a constructible function $\alp \in F(X)$,
$$c_*(\alp) = \chi(\alp) + \sum_i a_i[V_i]  \in H_*(X)$$
where $[V_i]$ is the homology class represented by a subvariety $V_i$ and $a_i \in \bZ$. Then the motivic class of the whole $c_*(\alp)$, denoted by $\mu c_*(\alp)$, could be defined by the following tautological one:
$$\mu c_*(\alp) = \chi(\alp)\jeden  + \sum_i a_i[V_i]  \in K_0(\Cal V_{\bC}).$$
However, this definition is not well defined. So it remains to see whether or there exists a motivic version $\mu c_*: F(X) \to K_0(\Cal V_{\bC})$. If we can obtain a reasonable motivic version $\mu c_*$, then we expect that we can also get a proalgebraic version of $\mu c_*$.
\endremark

\head \S 6 A few remarks on integrations \endhead

In this section we discuss a proalgebraic version of integration with respect to the proalgebraic Euler--Poincar\'e characteristic $\chi^{\op {pro}}$ or the Grothendieck ``motivic" class $\Ga^{\op {pro}}$ in the case when the directed set $\La = \bN$.

As remarked in the Introduction, the usual Euler--Poincar\'e characteristic $\chi: F(X) \to \bZ$
is described by, putting emphasis on integration, 
$$\chi(\alp) = \int_X \alp d\chi = \sum_n n \chi \bigl (\alp^{-1}(n) \bigr ).\tag 6.1$$
and furthermore, for a function $f: \bZ \to \bZ$ 
$$\int_X f(\alp) d\chi = \sum_n f(n) \chi \bigl (\alp^{-1}(n) \bigr ). \tag 6.2 $$
We remark that $f(\alp)$ should be more precisely expressed as the composite $f \circ \alp$, but that we write it so for simplicity .

A very na\"\i ve proalgebraic version of (6.1) would be
$$\int_{X_{\infty}} \alp d\chi^{\op {pro}} = \sum_n n \chi^{\op {pro}} \bigl (\alp^{-1}(n) \bigr )$$
for a proconstructible function $\alp \in F^{\op {pro}}(X_{\infty})$. 
Here comes out a problem. How do we define $\chi^{\op {pro}} \bigl (\alp^{-1}(n) \bigr )$ ? So far, $\chi^{\op {pro}}$ is defined on proconstructible functions or more explicitly defined on procharacteristic functions. As observed in Proposition (5.7), the procharacteristic function is equivalent to the cylinder set. Therefore $\chi^{\op {pro}} \bigl (\alp^{-1}(n) \bigr )$ are defined on cylinder sets (cf. Remark (5.9)).  However, we know that an arbitrary proconstructible function $\alp \in F^{\op {pro}}(X_{\infty})$ does not necessarily satisfy the property that $\alp^{-1}(n)$ is always a cylinder set for any $n$. For example, consider the following situation:
Let $\{X_n, \pi_{n(n+1)}:X_{n+1} \to X_n \}$ be a projective system of algebraic varieties with each structure morphism $\pi_{n(n+1)}:X_{n+1} \to X_n$ being surjective. Let $W_1 \in X_1$ and $W_2 \in X_2$ be constructible sets and consider the proconstructible function
$\alp = [\jeden_{W_1}] + [\jeden_{W_2}]$. Then we have
$$
\alp^{-1}(n) = \cases
\pi_1^{-1}(W_1) \ominus \pi_2^{-1}(W_2) \quad \text {if $n = 1$} \\
\pi_1^{-1}(W_1) \cap \pi_2^{-1}(W_2) \quad \text {if $n = 2$} \\
\emptyset \quad \text {otherwise} \endcases
$$
Here the symbol $\ominus$ is the symmetric difference, i.e., $A \ominus B = (A \setminus B) \cup (B \setminus A) = (A \cup B) \setminus (A \cap B).$

If $\pi_1^{-1}(W_1) \cap \pi_2^{-1}(W_2) = \emptyset$, then $\pi_1^{-1}(W_1) \ominus \pi_2^{-1}(W_2)$ is the disjoint union of two cylinder sets. Thus we can define $\chi^{\op {pro}} \bigl (\alp^{-1}(1) \bigr )$. However, if $\pi_1^{-1}(W_1) \cap \pi_2^{-1}(W_2) \not = \emptyset$, in general $\pi_1^{-1}(W_1) \ominus \pi_2^{-1}(W_2)$ and $\pi_1^{-1}(W_1) \cap \pi_2^{-1}(W_2)$ are not necessarily cylinder sets. Hence, to define the integration with respect to $\chi^{\op {pro}}$ or $\Ga^{\op {pro}}$, we need to restrict ourselves to proconstructible functions whose fibers are always cylinder sets. 

For a proalgebraic variety $\displaystyle X_{\infty} = \varprojlim _{n \in \bN} X_n$ and a projective system $\{p_{n m} \}$ of non-zero integers, consider a $\chi$-stable proconstructible
function $\alp$ whose fibers are all cylinder sets and a function $f: \bZ \to \bQ$. Then the following integration is well-defined:

$$\int_{X_{\infty}}f(\alp) d\chi^{\op {st.pro}}_{\{p_{n m} \}}:= \sum_k f(k) \chi^{\op {st.pro}}_{\{p_{n m} \}}\bigl (\alp^{-1}(k)  \bigr ).$$

Similarly, for a proalgebraic variety $\displaystyle X_{\infty}= \varprojlim _{n \in \bN} X_n$ and a projective system $\bigl \{[F_{nm}] \bigr \}$ of Grothendieck classes, consider a $\Ga$-stable proconstructible function $\alp$ whose fibers are all cylinder sets and a function $f: \bZ \to K_0(\Cal V_{\bC})_\Cal F$. Then the following ``motivic" integration is well-defined:

$$\int_{X_{\infty}}f(\alp) d\Ga^{\op {st.pro}}_{\{[F_{nm}] \}}:= \sum_k f(k) \Ga^{\op {st.pro}}_{\{[F_{nm}] \}}\bigl (\alp^{-1}(k)  \bigr ).$$

\head \S 7 proresolutions \endhead

Finally we consider the following projective system of resolution of singualrities:
Let $Y$ be a possibly singular variety and let $\Cal{RES}_Y$ be the collection $\{(Y', g) \}$ of resolution of singularities $g: Y' \to Y$ of Y, where $Y'$ is nonsingular and $g|_{Y' \setminus g^{-1}(Y_{\op {sing}})} : Y' \setminus g^{-1}(Y_{\op {sing}}) \to Y \setminus Y_{\op {sing}}$ is an isomorphism with $Y_{\op {sing}}$ denoting the singular set of $Y$. When $Y$ is nonsingular, $\Cal{RES}_Y$ is defined to be  just $\{(Y, \op {id}_Y) \}$. The reason for this requirement is that otherwise $\Cal{RES}_Y$ consists of all automorphisms of $Y$, which 
is not necessary for our purpose.

For two elements $g_1: Y_1 \to Y$, $g_2: Y_2 \to Y$ of $\Cal{RES}_Y$  we define the order $\leq$
as follows:  $g_1 \leq g_2$  if and only if there exists a morphism $g_{12}: Y_2 \to Y_1$ such that $g_2 = g_1 \circ g_{12}:$

$$\xymatrix{
Y_2 \ar[dr]_ {g_2}\ar[rr]^ {g_{12}} && Y_1\ar[dl]^{g_1}\\
& Y
}$$

\proclaim {Proposition (7.1)} For a possibly singular variety $Y$, the ordered set $(\Cal{RES}_Y, \leq)$ is a directed set.
\endproclaim

\demo {Proof}Let $g_1, g_2 \in \Cal{RES}_Y$. We want to show that there exists a desingularization
$g_3:Y_3 \to Y$ such that $g_1 \leq g_3$ and $g_2 \leq g_3$. Consider the fiber product

$$\CD
Y_1 \times _Y Y_2 @> {\widetilde {g_1}} >> Y_2 \\
@V {\widetilde {g_2}}VV @VV {g_2}V\\
Y_1@>> {g_1} > Y, \endCD
$$

and furthermore consider a resolution of singularities of $Y_1 \times _Y Y_2$:

$$\pi : \widetilde {Y_1 \times _Y Y_2} \to Y_1 \times _Y Y_2.$$
Let us set $Y_3$ to be $\widetilde {Y_1 \times _Y Y_2}$, and we set
$g_3: Y_3 \to Y$ to be the composite 

$$g_3 := g_1 \circ \widetilde {g_2} \circ \pi = g_2 \circ \widetilde {g_1} \circ \pi.$$
Which means that $g_1 \leq g_3$ and $g_2 \leq g_3$.
\qed
\enddemo

Therefore we can see that the collection $\Cal{RES}_Y$  of resolution of singularities becomes a projective system with the index set being the directed set $\Cal{RES}_Y$ itself.

\definition {Definition (7.2)} For a possbily singular variety $Y$, the projective limit of the 
projective system $\Cal{RES}_Y$ of resolutions of singularities of $Y$
is called {\it the proresolution of $Y$}(a sort of ``maximal" resolution of X) and denoted by
$$\pi^{\op {pro}}: \widetilde Y^{\op {pro}} \to Y.$$
\enddefinition

This proresolution is motivated by our previous work [Y4].
Our na\" \i ve question is a relationship between the Nash arc space and the proresolution:

\proclaim {Question (7.3)} Let $X$ be a possibly singular variety. Then does there exist a canonical promorphsim $\nu: \Cal L(X) \to \widetilde X^{\op {pro}}$ as a $X$-provariety ?, i.e, such that the following diagram commutes

$$\xymatrix{
{\Cal L(X)}  \ar[dr]_{} \ar[rr]^ {\nu} && {\widetilde X^{\op {pro}}} \ar [dl]_{}\\
& X.
}$$
\endproclaim

One could ask if there exists a canonical morphism from $\widetilde X^{\op {pro}}$ to $\Cal L(X)$
over $X$, but it would be reasonable to consider the above question, because over the nonsingular part of $X$ the proresolution $\widetilde X^{\op {pro}}$ is isomorphic to the nonsingular part of $X$, whereas the Nash arc space $\Cal L(X)$ is a fiber bundle over the nonsinguar part $X_{\op {smooth}}$ of $X$ with fiber being the proalgebraic variety $\bC^{\bN}$, to be more precise, 
$X_{\op {smooth}}\widetilde {\times} (\bC^{\op {dim}X})^{\bN}$, using the notation given in \S 5.  We hope to return to this question.

\enddocument